\pdfoutput=1
\documentclass[ijoc,nonblindrev]{informs3} 

\DoubleSpacedXII

\usepackage{ulem} 
\usepackage{bbm} 
\usepackage{mathtools} 
\usepackage{xcolor} 
\usepackage{algorithm,algcompatible} 
\usepackage[labelformat=simple, size=scriptsize]{subcaption} 

\usepackage{multirow} 
\usepackage{booktabs} 
\usepackage{hyperref} 
\usepackage{amsfonts} 

\usepackage{natbib}
 \bibpunct[, ]{(}{)}{,}{a}{}{,}%
 \def\bibfont{\small}%
\TheoremsNumberedThrough     

\EquationsNumberedBySection 

\begin{document}


\RUNAUTHOR{Baik, Byon, and Ko}

\RUNTITLE{DR Stratified Sampling under Multiple Input Uncertainties}

\TITLE{Distributionally Robust Stratified Sampling for Stochastic Simulations with Multiple Uncertain Input Models}

\ARTICLEAUTHORS{%
\AUTHOR{Seung Min Baik}
\AFF{Pohang University of Science and Technology, \EMAIL{seungmin.baik@lstlab.org}}
\AUTHOR{Eunshin Byon}
\AFF{University of Michigan, \EMAIL{ebyon@umich.edu}, \URL{https://ebyon.engin.umich.edu/}}
\AUTHOR{Young Myoung Ko}
\AFF{Pohang University of Science and Technology, \EMAIL{youngko@postech.ac.kr}, \URL{https://www.lstlab.org/}}
} 

\ABSTRACT{
This paper presents a robust version of the stratified sampling method when multiple uncertain input models are considered for stochastic simulation. Various variance reduction techniques have demonstrated their superior performance in accelerating simulation processes. Nevertheless, they often use a single input model and further assume that the input model is exactly known and fixed. We consider more general cases in which it is necessary to assess a simulation's response to a variety of input models, such as when evaluating the reliability of wind turbines under nonstationary wind conditions or the operation of a service system when the distribution of customer inter-arrival time is heterogeneous at different times. Moreover, the estimation variance may be considerably impacted by uncertainty in input models. To address such nonstationary and uncertain input models, we offer a \textit{distributionally robust (DR) stratified sampling} approach with the goal of minimizing the maximum of worst-case estimator variances among plausible but uncertain input models. Specifically, we devise a bi-level optimization framework for formulating DR stochastic problems with different ambiguity set designs, based on the $L_2$-norm, 1-Wasserstein distance, parametric family of distributions, and distribution moments. 
In order to cope with the non-convexity of objective function, we present a solution approach that uses Bayesian optimization.
Numerical experiments and the wind turbine case study demonstrate the robustness of the proposed approach. 

}%



\maketitle

%



\section{Introduction}
\label{sec:introduction}

This paper devises a new input sampling strategy for stochastic simulation to estimate outputs of interests under multiple uncertain input models. To acquire system outputs, stochastic simulation typically samples random input parameters and runs a computer model repeatedly. We focus on stochastic computer models that produce noisy outputs despite identical input parameters~\citep{choe2015importance}. Stochastic simulation with stochastic computer models involve two sources of randomness: the input's probability distribution and the output's inherent stochasticity.

This study is motivated from reliability analysis for wind turbines using stochastic simulation~\citep{Choe2016}. The National Renewable Energy Laboratory (NREL) of the U.S. Department of Energy has created aeroelastic computer models, such as TurbSim~\citep{jonkman2009turbsim} and FAST~\citep{jonkman2005fast}, to aid in the design of reliable wind turbines. To analyze the failure probability, that the load response exceeds a certain threshold level, variance reduction techniques have been proposed to enhance computing efficiency over the crude Monte Carlo sampling~\citep{choe2015importance,pan2020adaptive,myoung2022optimal}. 

Variance reduction studies often employ a single input model. Yet, some situations require handling several input models, such as when the input characteristics change over time or across dispersed locations. Consider a multi-turbine wind farm. Each turbine experiences a different wind condition because upstream turbines' operations add to the turbulence, which changes the free-flow wind condition~\citep{You2016,liu2022parameter}, referred to as \textit{wake effects}~\citep{YLB2018}. As a result, downstream turbines experience heterogeneous wind conditions. Furthermore, even at a fixed location, the wind patterns change throughout the year~\citep{jang2020probabilistic}. Calculating the failure probabilities by running simulations with various input models will require extensive computing power. On the other hand, the optimal sampling budget allocation for a specific input model may result in significant inefficiency for other input models.

Furthermore, conventional variance reduction approaches assume that the \textit{true} input distribution is known. But occasionally, a fitted or empirical distribution that is estimated with limited observations is used as its surrogate. When measurement data is unavailable, a physics-based numerical model is employed to approximate the true distribution~\citep{zhang2015comparison}. The estimation errors in the input model may result in poor estimation quality of the simulation response. Though many studies have been carried out recently to take input uncertainty into account, the majority are yet limited to a single input model.

Among several variance reduction techniques, this study is concerned with stratified sampling. We devise a new variance reduction technique, referred to \textit{distributionally robust stratification} (shortly, \textit{DR-strat}), for determining a robust input sampling strategy. Our approach involves allocating the limited simulation budget when estimating performance measures under different uncertain input models. Hinging upon the fundamentals of distributionally robust optimization (DRO), we minimize the worst-case estimator variance for a set of plausible distributions. Specifically, we formulate a bi-level optimization problem where the outer problem minimizes the maximum variance using a sampling vector across strata as a decision vector, while the inner problem finds a plausible (uncertain) input model with the largest variance. We employ Bayesian optimization (BO) to search the solution space probabilistically.

Below we summarize the contribution of our study. 
\begin{itemize}	
	\item We propose a new variance reduction technique to determine a robust input sampling strategy under multiple input models' uncertainties. To the best of our knowledge, this is the first study to take the input model uncertainty into account in variance reduction techniques for stochastic simulation.

    \item We provide a framework for formulating an optimization problem to derive a robust input sampling strategy. In contrast to most existing DRO studies, which deal with a single input model, we consider multiple input models in formulating the bi-level DR stochastic problem and suggest a solution procedure by adopting BO.

	\item We construct four types of ambiguity sets of plausible distributions that represent potential candidates for true input models, based on $L_2$-norm, 1-Wasserstein distance, parametric family, and distribution moments. We also investigate how various set design approaches affect the estimation result.  While we demonstrate four types of ambiguity sets, the proposed bi-level optimization methodology is easily extensible to other types of sets as well.
	
    \item Our numerical experiments and a case study involving wind turbine reliability demonstrate that the proposed method successfully derives an estimator that robustly reacts to various input model uncertainties. As a result, our approach enables the efficient reuse of simulation results for performance measure estimation under multiple uncertain input models, which is crucial in the circumstances with limited computational budgets.
\end{itemize}  

The remainder of the paper is organized as follows. Section~\ref{sec:literature_review} reviews previous studies. 
Section~\ref{sec:problem_description} summarizes the conventional stratification method and provides the overall framework for deciding a robust input sampling strategy. Section~\ref{sec:methodology} discusses the DR-stratified sampling method with mathematical details. Section~\ref{sec:numerical_experiments} conducts numerical experiments. Section~\ref{sec:conclusions} concludes and suggests future research directions.

\section{Literature Review}
\label{sec:literature_review}

Overall, this study is closely related to the two broad areas of research: stochastic simulation under input uncertainty and DRO. First, studies on input uncertainty in stochastic simulation include multiple research streams, including input uncertainty quantification~\citep{song2017input}, sensitivity analysis on additional input data collection, and the simulation optimization under input uncertainty~\citep{zhou2017simulation}. \cite{CORLU2020100162} and \cite{Barton2022} provide a comprehensive review of related research studies.

This work is more closely related to the third of these streams. In particular, our approach is similar to the computational budget allocation problem in ranking and selection (R\&S) studies, concerning input uncertainty to pursue a robust optimal sampling method. \cite{song2015input} investigate the impact of input uncertainty on simulation output with a mixed-effect model and adjust indifference-zone (IZ) procedures to guarantee the average probability of correct selection (PCS). \cite{gao2017new} follow a robust approach for optimal computing budget allocation and solve approximate optimization problems to maximize PCS. \cite{fan2020distributionally} study a robust selection of the best problem with the IZ approach based on the concept of an ambiguity set. 

When evaluating performance, these approaches consider both the alternative and the input model to search for the best among a set of alternatives. Our focus is slightly different, as we are particularly interested in the performance (i.e., estimator variance) solely impacted by the input model. While R\&S studies typically require separate simulations for different alternatives under the same input model, which is both effective and necessary for their purposes, we simultaneously assess the influence of the sampling strategy across multiple input models. 

Next, studies in the DRO literature treat uncertain input models with the concept of ambiguity set. \cite{delage2010distributionally} conduct a study on modeling the DRO problem with a moment-constrained ambiguity set and developing a tractable solution procedure for solving it. \cite{lam2016empirical} adopt the empirical likelihood method to interpret the conventional DRO approach and investigate the confidence interval for the target performance to address a potential loss of coverage accuracy. \cite{lam2017tail} estimate the tail-related quantity of interest and investigate the characteristics of the worst-case objective. \cite{rahimian2019distributionally} review related studies comprehensively. 

Similar to these DRO studies that generate an ambiguity set, our approach makes use of ambiguity sets. However, we consider several ambiguity sets, one corresponding to each input model, unlike most previous DRO studies that only analyze a single input model.


\section{Problem Description}
\label{sec:problem_description}

Consider a black box computer model that generates an output $Y \in \mathbb{R}$ given an input $X \in \mathbb{R}^{P}$ following a distribution \(F\). Given \(X\), the computer model produces either a stochastic (or noisy) or deterministic output. In this study, we focus on the stochastic computer model, mirroring the stochasticity of NREL simulators employed in our motivating wind turbine application. However, our approach can be easily adopted in deterministic computer models. 

Let \(Y(X)\) denote the simulation output at the input  \(X\). For the reliability analysis to estimate a failure probability \(\mathbb{P}(Y(X)>l)\), representing the probability of the simulation output being larger than a threshold $l$, we use \(g(x) = \mathbbm{1}(Y(x)>l)\) where \(\mathbbm{1}(\cdot)\) is an indicator function. With the stochastic computer model, \(g(x)\) is random even at fixed \(x\). We are interested in estimating the mean of \(g(X)\) (i.e., \(\mu \coloneqq \mathbb{E}\left[ g(X)\right] = \mathbb{E}_X \left[ \mathbb{E}_Y \left[g(X) | X \right] \right]\)). With $g(x) = \mathbbm{1}\left(Y(x) > l\right)$, we have \(\mu = \mathbb{E}_X \left[ \mathbb{E}_Y \left[\mathbbm{1}\left( Y(X)>l|X\right) \right] \right] = \mathbb{E}_X \left[ \mathbb{P} \left( Y(X)>l|X\right) \right] = \mathbb{P}\left(Y(X)>l\right)\). Our objective is to design an estimator $\hat{\mu}$ that effectively estimates the target performance measure $\mu$. Proper allocation of the simulation efforts is crucial under a fixed computational budget \(N_T\) when the computational cost for evaluating \(g(X)\) is expensive. 

This study considers a discrete input vector \(X\), as a starting point of research that addresses input uncertainty in variance reduction techniques, for computational purpose. Most existing DRO studies has primarily focused on ensuring tractability when constructing ambiguity sets and formulating optimization problems~\citep{rahimian2019distributionally}. Considering a continuous input has often led to situations where the optimization problem becomes computationally intractable, except for special cases with inherent structural features. Therefore, many practical problems have assumed a discrete input~\citep{bansal2018}, as it allows for feasible solution procedure to DRO problem. Similarly, we also use a discrete input so to ensure that our DR-strat problem can be solved under all four types of ambiguity sets.

Still, we would like to note that the proposed methodology is practically applicable to situations where discretization of continuous inputs can be employed. One of the most commonly used methods in the literature on wind energy reliability is the so-called \textit{binning method}~\citep{Choe2016}. It partitions the wind speed range into multiple intervals and runs a computer model at each interval (or bin). The strata in the stratified sampling can be formed by these intervals, and their representative values can be set to be the domain of $X$. Further, when the input vector is continuous (e.g., wind speed), we can discretize it into multiple bins, as demonstrated in our case study in Section~\ref{sec:numerical_experiments}.

\subsection{Recap: Stratified Sampling for Single Input Model}
\label{subsec:recap_stratification}

The crude Monte Carlo sampling is the most basic approach that provides an unbiased estimator for a single input model with the distribution $F$. It estimates the performance measure \({\mu}\) by \(\hat{\mu}^{\text{MC}} = \sum_{n=1}^{N_T} g(X_n)/N_T\), where $\{X_n\}_{n = 1}^{N_T}$ are independent and identically distributed (i.i.d.) samples drawn from $F$. When the event of interest occurs rarely, such as the exceedance event \( \{Y(x)>l\} \) with large \(l\), a significant number of simulation runs may be required to obtain an accurate output estimate. Alternatively, stratified sampling, one of the popular variance reduction techniques, provides more effective way for drawing input samples to reduce the estimator variance $\text{Var}[\hat{\mu}]$.

Let us review conventional stratified sampling for a single input model.  Suppose that the sampling domain \(\Omega\) of the input vector \(X\) can be divided into  mutually exclusive and exhaustive strata \(\{S_k\}_{k= 1}^{K}\). Conditional output mean \(m_k = \mathbb{E}\left[g(X)|X\in  S_k\right]\) of the $k$th stratum can be estimated by averaging the simulation outputs at \(n_k\) conditional inputs as \(\sum_{j=1}^{n_k} g\left(X_{j|k}\right)/n_k\), where \(\{X_{j|k}\}_{j = 1}^{n_k}\) are i.i.d. samples drawn from a conditional distribution of $F$ given that an input belongs to the $k$th stratum (i.e., $\{X\in S_k\}$). We call $\boldsymbol{n} = (n_1, n_2, \dots, n_K)$ a sampling vector. A probability of the $k$th stratum is $\omega_k = \mathbb{P}(X \in S_k)$. We assume $\omega_k > 0, \forall k$, to avoid any trivial issues. With strata probabilities $\boldsymbol{\omega} = (\omega_1, \omega_2, \dots, \omega_K)$, we define the stratified sampling estimator by aggregating the conditional estimates from each of the strata as follows:
\begin{equation}
    \hat{\mu}^{\text{Str}}(\boldsymbol{n}) = \sum_{k=1}^K \omega_k \frac{\sum_{j=1}^{n_k} g\left(X_{j|k}\right)}{n_k}.
    \label{eq:estimator_strat}
\end{equation}		

The stratification estimator is always unbiased (i.e., $\mathbb{E} [\hat{\mu}^{\text{Str}}(\boldsymbol{n})]=\mu$) regardless of $\boldsymbol{n}$. However, the sample vector $\boldsymbol{n}$ affects the stratified sampling estimator variance. Please refer Online Supplement A.1 for details. Suppose that the computational cost of drawing an input, as well as evaluating an output, is the same across all strata. Given a total simulation budget \(N_T\), the following sampling vector is known to be optimal (i.e., it minimizes $\text{Var}[\hat{\mu}^{\text{Str}}(\boldsymbol{n})]$)~\citep{owen2013monte}:
\begin{equation*}
    \boldsymbol{n}^{\text{Str}} = (n_1^{\text{Str}}, n_2^{\text{Str}}, \dots, n_K^{\text{Str}}), \ n_k^{\text{Str}} = N_T \frac{\omega_k \sigma_k}{\sum_{k=1}^K \omega_k \sigma_k}, \ \forall k = 1, 2, \dots, K.
\end{equation*}
In practice, $n_k^{\text{Str}}$'s are rounded to integers by allowing small non-proportionalities.

\subsection{Multiple Input Models}
\label{subsec:non_stationary_input_model}

As discussed earlier in Section~\ref{sec:introduction}, multiple uncertain input models need to be taken into account in several circumstances. Considering $M$ input models, we are interested in estimating $M$ performance measures $\mu_m^c \coloneqq \mathbb{E}\left[ g(X_m^c) \right]$ for $m=1, 2, \dots, M$, where the input random vector or variable (R.V.) $X_m^c$ follows the $m$th input distribution $F_m^c$. Here, the superscript $c$ is used to denote the \textit{correct} (or \textit{true}) information. We suppose that these true distributions are unknown, and they are inferred using empirical data in practice. Our goal is to design an estimator that performs well for all $M$ input models in terms of reducing variance while also being robust to the uncertainties in input distributions.

Specifically, we minimize the maximum of $M$ estimator variances, $\underset{1 \leq m \leq M}{\max} \text{Var}\left[\hat{\mu}_m (\boldsymbol{n}) \right]$, where each variance corresponds to the estimator for a different input model with the same sampling vector $\boldsymbol{n}$. However, the precise maximum value is impossible to calculate because the input distributions are not known. To tackle this, we consider an ambiguity set, denoted by $\mathcal{F}_m$, to represent the set of probable distributions of the $m$th input model for $1 \leq m \leq M$.  A set $\mathcal{F}_m$ is constructed to include distributions \textit{close} to a nominal distribution (e.g., predicted or fitted input distribution), such as those within a certain Wasserstein distance. Consequently, to take the robustness against the uncertainty, we adopt a DRO approach, whereby we consider the worst-case estimator variance over a set of distributions.

We aim to allocate computational budgets across strata to minimize the maximum variance among multiple uncertain input models. Let $\hat{\mu}^{\text{DR-Str}}(\boldsymbol{n};F_m)$ denote the estimator for the $m$th input model under our robust stratified sampling approach (the mathematical definition of $\hat{\mu}^{\text{DR-Str}}(\boldsymbol{n};F_m)$ will be provided in Section~\ref{sec:methodology}). Then the problem boils down to finding a robust input sampling vector $\boldsymbol{n}^{\text{DR-Str}}$, that is, how many samples to draw from the predetermined strata, by solving the following problem.
\begin{equation}
    \underset{\boldsymbol{n}}{\min} \ \underset{1 \leq m \leq M}{\max} \ \underset{F_m \in \mathcal{F}_m}{\max} \   \text{Var} \left[\hat{\mu}^{\text{DR-Str}}(\boldsymbol{n};F_m) \right].
    \label{eq:minmax_original}
\end{equation}
By minimizing the maximum of worst-case estimator variances, we prevent the estimator variance from growing too large even with the poor estimation of uncertain input models.

\section{Methodology: DR-Stratified Sampling}
\label{sec:methodology}

The conventional stratification estimator, discussed in Section~\ref{subsec:recap_stratification}, is determined based on the characteristics of input distribution \(F\) and output function \(g\). Thus, the optimal simulation budget allocation, or the sampling vector $\boldsymbol{n}^{\text{Str}}$, changes as \(F\) varies. 
This section proposes a distributionally robust stratification method designed to robustly respond to uncertainties within multiple input models. 

\subsection{Formulation of DR-Strat Problem}
\label{subsec:DR_stratification_formulation}

This section presents the detailed formulation of the DR-strat problem for determining the DR-strat sampling vector in \eqref{eq:minmax_original}. We first define the new estimator design that is suitable to handle multiple uncertain input models. Then we formulate a bi-level optimization problem where the inner problem finds the worst-case estimator variance among the plausible input models and the outer finds the optimal sampling vector.

\subsubsection{DR-Stratfication Estimator.}
\label{subsubsec:DR_stratification_estimator}

To estimate outputs of interest under several input models, our strategy is to run simulations under a \textit{reference} distribution (a single common distribution used to draw inputs for all models), instead of running simulations under each input model (regarding individual sampling distribution) separately. We then reuse the obtained simulation outcomes for each input model. This procedure enables us to significantly reduce simulation efforts. The problem is how to allocate simulation efforts.
 
Let us consider the $m$th input model. The new estimator considers both the reference distribution $F_{ref}$ and the plausible distribution $F_m$, a candidate for characterizing the input model. Here, $F_m$ is an element of \(\mathcal{F}_m\) constructed upon the nominal (or base) distribution $\bar{F}_m$ of the $m$th input model. This nominal distribution serves as a basis (such as the center point) for creating the ambiguity set. For the input sampling domain $\Omega = \{x_i\}_{i=1}^{|\Omega|}$, we denote the probability mass function (pmf) values of $F_{ref}$ as $\boldsymbol{p}_{ref} = \left( p_{ref,1}, \dots, p_{ref, |\Omega|} \right)$ and $F_m$ as $\boldsymbol{p}_m = \left( p_{m,1}, \dots, p_{m,|\Omega|} \right)$. Further, we use the notations $X_{ref}$ and $X_m$ to denote input R.V.s following $F_{ref}$ and $F_m$, respectively. Thus, we have $\mathbb{P}\left( X_{ref} = x_i \right) = p_{ref, i}$ and $\mathbb{P}\left(X_m = x_i \right) = p_{m, i}$ for $1 \leq i \leq |\Omega|$. Suppose that we divide $\Omega$ into $K$ strata for $K \le |\Omega|$. The probabilities that  $X_{ref}$ and $X_m$ belong to the $k$th stratum, $S_k$, become $\omega_{ref,k} = \mathbb{P}\left( X_{ref} \in S_k \right) = \sum_{i \in \{i| x_i \in S_k\}} p_{ref,i}$ and $\omega_{m,k} = \mathbb{P}\left( X_m \in S_k \right) = \sum_{i \in \{i| x_i \in S_k\}} p_{m,i}$ for $k=1,2,\dots,K$. Similar to the conventional stratified sampling, we assume these strata probabilities are strictly positive to avoid trivial issues. 

We additionally define notations $X_{ref,k}$ and $X_{m,k}$ to denote the conditional R.V.s, given that $X_{ref}$ and $X_m$ belong to the $k$th stratum, respectively (i.e., $X_{ref,k} \buildrel d \over = X_{ref}|\{X_{ref} \in S_k\}$ and $X_{m,k} \buildrel d \over = X_m|\{X_m \in S_k\}$). So, the conditional probabilities of input $x_i$ given that $\{x_i \in S_k\}$ are $\mathbb{P} \left( X_{ref,k} = x_i \right) = p_{ref,i} / \omega_{ref,k}$ and $\mathbb{P}\left( X_{m,k} = x_i \right) = p_{m,i} / \omega_{m,k}$.

Suppose we draw $n_k$ i.i.d. samples, denoted by  $\{X_{j|k}\}_{j=1}^{n_k}$, from the conditional reference distribution of $F_{ref}$ given $\{ X_{ref} \in S_k \}$ for $1 \leq k \leq K$. Then, the new DR-strat estimator for estimating $E[g(X_m)]$ under the  input distribution $F_m$ can be defined as follows:
\begin{equation}
    \hat{\mu}^{\text{DR-Str}} (\boldsymbol{n}; F_m) = \sum_{k=1}^K \frac{\omega_{m,k}}{n_k} \sum_{j=1}^{n_k} g\left(X_{j|k}\right) \frac{\mathbb{P}\left(X_{m,k} = X_{j|k} \right)}{\mathbb{P}\left(X_{ref,k} = X_{j|k} \right)} .
    \label{eq:dr_strat_estimator}
\end{equation}
Here, please note that $\hat{\mu}^{\text{DR-Str}} (\cdot)$ has an additional argument $F_m$ (for the evaluation), unlike $\mu^{\text{Str}}(\cdot)$ in \eqref{eq:estimator_strat} that does not. This indicates that the estimation is performed for the $m$th input model with distribution $F_m$. Further, the strata probability $\omega_k$ in \eqref{eq:estimator_strat} is substituted with $\omega_{m,k}$ in order to consider $F_m$ in the left-hand side of the equation. 

We would like to highlight that there is another important difference between the estimators $\hat{\mu}^{\text{Str}}$ and $\hat{\mu}^{\text{DR-Str}}$. In estimating the measure of interest when the same distribution is used for both input sampling and evaluation, $\hat{\mu}^{\text{Str}}$ in~\eqref{eq:estimator_strat} provides unbiased estimation for $E[g(X)]$. On the contrary, DR-strat samples inputs from the reference distribution $F_{ref}$ but estimates the output under another distribution $F_m$. Thus, we need to use the likelihood ratio  $\left. \mathbb{P}\left(X_{m,k} = X_{j|k} \right) \middle/ \mathbb{P}\left(X_{ref,k} = X_{j|k}\right) \right.$ in $\hat{\mu}^{\text{DR-Str}}$ in~\eqref{eq:dr_strat_estimator} to correct the bias.

Proposition~\ref{eq:proposition_estimator_mean} shows that the DR-strat estimator is unbiased (i.e., the estimator mean becomes the same as the true output mean when the input distribution is $F_m$).
\begin{proposition}
    \label{eq:proposition_estimator_mean}
    For a random vector $X_m$ following a distribution $F_m$,
    \begin{equation}
        \mathbb{E}\left[\hat{\mu}^{\text{DR-Str}} (\boldsymbol{n};F_m) \right] = \mathbb{E} \left[ g(X_m) \right].
        \label{eq:unbiasedness_of_dr_strat_estimator}
    \end{equation}
\end{proposition} 
Next, Proposition~\ref{eq:proposition_estimator_variance} derives the variance of the DR-strat estimator.
\begin{proposition}
    \label{eq:proposition_estimator_variance}
    For a random vector $X_m$ following a distribution $F_m$,
    \begin{equation}
        \begin{aligned} 
            &\text{Var} \left[ \hat{\mu}^{\text{DR-Str}} (\boldsymbol{n};F_m) \right] \\
            &=  \sum_{k=1}^{K} \frac{1}{n_k} \left( \left( \omega_{ref,k} \sum_{i \in \{i| x_i \in S_k\}} \mathbb{E}\left[ g(x_i) \right] \frac{ { \mathbb{P}\left(X_{m} = x_i \right) }^2}{ \mathbb{P}\left(X_{ref} = x_i \right) } \right) - \left( \sum_{i \in \{i| x_i \in S_k\}} \mathbb{E}\left[ g(x_i) \right] \mathbb{P}\left(X_{m} = x_i \right) \right)^2 \right).
        \end{aligned}
        \label{eq:dr_strat_estimator_variance}
    \end{equation}
\end{proposition} 
Online Supplement A.2 and A.3 provide the detailed proofs for the above propositions.

As $\hat{\mu}^{\text{DR-Str}} (\boldsymbol{n};F_m)$ is an unbiased estimator for $\mathbb{E} \left[ g(X_m) \right]$ as shown in \eqref{eq:unbiasedness_of_dr_strat_estimator}, we want to minimize its variance $\text{Var} \left[ \hat{\mu}^{\text{DR-Str}} (\boldsymbol{n};F_m) \right]$ in \eqref{eq:dr_strat_estimator_variance} by allocating simulation budgets adequately. In the subsequent discussion, we will present a new formulation to robustly allocate budgets across multiple strata in order to handle multiple uncertain input distributions.

\subsubsection{DR-Strat Problem.}
\label{subsubsec:DR-Stratification_problem}

We start constructing the DR-strat problem by formulating the inner maximization problem first, which aims to find the maximum value of worst-case variances among multiple sets of plausible input models. In this stage, the sampling vector $\boldsymbol{n}$ (the decision vector of the outer minimization problem) is given. Using the new estimator design in~\eqref{eq:dr_strat_estimator}, the inner maximization problem becomes
\begin{equation*}
    \underset{1 \leq m \leq M}{\max} \ \underset{F_m \in \mathcal{F}_m}{\max} \ \text{Var} \left[\hat{\mu}^{\text{DR-Str}}(\boldsymbol{n};F_m) \right], 
\end{equation*}
where the plausible distributions $F_m$'s and input model index $m$ are the decision variables.

We define an index set of input values at the $k$th stratum as $I_k = \{i | x_i \in S_k\}$ for $1 \leq k \leq K$. Noting that $X_m$ is a discrete R.V., the distribution $F_m \in \mathcal{F}_m$ has the equivalent meaning with $\boldsymbol{p}_m \in \mathcal{P}_m$ with  $\mathcal{P}_m$ being the ambiguity set expressed in terms of pmfs. Using the estimator variance in \eqref{eq:dr_strat_estimator_variance}, the inner problem can be reformulated as follows:
\begin{equation}
    \underset{1 \leq m \leq M}{\max} \ \underset{ \boldsymbol{p}_m \in \mathcal{P}_m}{\max} \ \sum_{k=1}^{K} \frac{1}{n_k} \left( \left( \omega_{ref,k} \sum_{ i \in I_k } \mathbb{E}\left[ g(x_i) \right] \frac{ p_{m,i}^2}{ p_{ref,i} } \right) - \left( \sum_{ i \in I_k} \mathbb{E}\left[ g(x_i) \right] p_{m,i} \right)^2 \right).
   \label{eq:inner_problem_pmf_version}
\end{equation}
Here, $\mathbb{E}[g(x_i)]$ in \eqref{eq:inner_problem_pmf_version} are supposed to be estimated from the pilot stage simulation (e.g., by fitting meta-models to data).

Next, to find the optimal sampling strategy that minimizes the maximum value of worst-case estimator variances, we formulate the DR-strat problem as follows:
\begin{equation}
    \begin{aligned}
        \text{(DR-Str)} \ \underset{\boldsymbol{n}}{\min}
        \underset{\substack{\boldsymbol{p}_m \in \mathcal{P}_m \\ 1 \leq m \leq M}}{\max} 
        \ & \sum_{k=1}^{K} \frac{1}{n_k} \left( \left( \omega_{ref,k} \sum_{ i \in I_k } \mathbb{E}\left[ g(x_i) \right] \frac{ p_{m,i}^2}{ p_{ref,i} } \right) - \left( \sum_{ i \in I_k} \mathbb{E}\left[ g(x_i) \right] p_{m,i} \right)^2 \right)\\
      s.t. \ \ & \sum_{k=1}^K n_k = N_T
    \end{aligned}.
    \label{eq:min_max_opt_model}
\end{equation}
The outer minimization problem determines a sampling vector $\boldsymbol{n}$ under a budget constraint. We call the optimal solution of this min-max problem a \textit{DR-strat sampling vector}, denoted by $\boldsymbol{n}^\text{DR-Str}$.  
Section~\ref{subsec:solving_DR-Strat_problem} describes how we solve this problem.

\subsection{Ambiguity Set Design}
\label{subsec:uncertainty_set_design}

This section discusses the design of the ambiguity set in the DR-strat problem. The configuration of the ambiguity set significantly affects the result from the DR-strat approach, as it determines the search space of the inner maximization problem. We explore four types of ambiguity sets, those often employed in the literature~\citep{rahimian2019distributionally}. The first two sets are based on discrepancy functionals associated with the $L_2$-norm and 1-Wasserstein distance. We also construct an ambiguity set based on distribution moments. Finally, a collection of the same parametric distributions is employed. These four set types provide a comprehensive analysis of the DR-strat's performance under various aspects of input model uncertainty, while other set design can also be used, based on the specific problem structure at hand and the prior knowledge available about the input distribution.

We define an ambiguity set as a collection of pmfs, as we consider discrete input vectors $X_m$'s. Let $\bar{\boldsymbol{p}}_m = (\bar{p}_{m,1}, \bar{p}_{m,1}, \dots, \bar{p}_{m,|\Omega|})$ denote the $m$th nominal distribution. The elements in $\bar{\boldsymbol{p}}_m$ or $\boldsymbol{p}_m$ should add up to one (i.e., $\sum_{i=1}^{|\Omega|} \bar{p}_{m,i} = \sum_{i=1}^{|\Omega|} p_{m, i}$ = 1). We let a positive scalar value $\gamma_m$ denote a parameter that quantifies the degree of uncertainty. Depending on the set design, we will use an extra subscript or superscript in the subsequent discussion. The size parameter $\gamma_m$ can be chosen using domain knowledge or the level of confidence about the nominal distribution. 

We assume that each realization of the $m$th input model, ${\boldsymbol{p}_m} \in \mathcal{P}_m$, is independent of the realization of another input model, ${\boldsymbol{p}_{m'}} \in \mathcal{P}_{m'}$, when $m \neq m'$. Thus, we construct the ambiguity set for each input model separately. Future extension of this research may address possible dependencies between input models.

Now, we discuss each type of ambiguity set for the $m$th input model. First, we define the ambiguity set based on the $L_2$-norm as follows:
\begin{equation}
    \mathcal{P}_m^{L_2} = \left\{\boldsymbol{p}_m \middle\vert \ \left\|\boldsymbol{p}_m - \bar{\boldsymbol{p}}_m \right\|_{2} \leq \gamma_m^{L_2} \right\} = \left\{ \boldsymbol{p}_m \middle\vert \ \sum_{i=1}^{|\Omega|} \left( p_{m,i}-\bar{p}_{m,i} \right)^2 \leq \left(\gamma_m^{L_2}\right)^2  \right\},
    \label{eq:uc_l2_1d}
\end{equation}
where $\|\cdot\|_{2}$ denotes the $L_2$-norm. This set consists of pmfs that have the $L_2$ distance to the nominal pmf $\bar{\boldsymbol{p}}_m$ smaller than the uncertainty level $\gamma_m^{L_2}$. 

Next, the $p$-Wasserstein ($p$-$\mathcal{W}$) distance-based ambiguity set is defined as follows:
\begin{equation*}
    \mathcal{P}_m^{p-\mathcal{W}} = \left\{ \boldsymbol{p}_m \middle\vert \ \exists q_{m, ij} \geq 0, \ \forall i,j = 1, \dots,  |\Omega|, \ s.t. \
    \begin{array}{l} 
        \sum\limits_{i=1}^{|\Omega|} \sum\limits_{j=1}^{|\Omega|} \| x_i - x_j\|^{p} q_{m, ij} \leq \left( \gamma_m^{p-\mathcal{W}} \right)^p\\
        \sum\limits_{j=1}^{|\Omega|} q_{m,ij} = p_{m,i}, \ \forall i = 1, \dots, |\Omega| \\
        \sum\limits_{i=1}^{|\Omega|} q_{m, ij} = \bar{p}_{m,j}, \ \forall j = 1, \dots,  |\Omega|
    \end{array}
    \right\},
\end{equation*}
where $\|\cdot\|$ is the basis norm used for $p$-$\mathcal{W}$ distance. This set consists of pmfs $\boldsymbol{p}_m$'s, in which the $p$-$\mathcal{W}$ distance to the nominal distribution $\bar{\boldsymbol{p}}_m$ smaller than the uncertainty level $\gamma_m^{p-\mathcal{W}}$. Please refer Online Supplement A.4 for details.

Several studies suggest various techniques for solving DRO problems regarding Wasserstein distance-related constraints~\citep{rahimian2019distributionally}. For illustrative purposes, we present a case where such constraints are relatively easy to handle. When the input R.V. is defined on one-dimensional space (i.e., when $p=1$), the $1$-$\mathcal{W}$ distance-based ambiguity set with the size parameter $\gamma_m^{1-\mathcal{W}}$ becomes 
\begin{equation}
    \mathcal{P}_m^{1-\mathcal{W}} = \left\{ \boldsymbol{p}_m \middle\vert \  \sum_{i=1}^{|\Omega|-1} \left( \left| \sum_{j=1}^i p_{m,j} - \sum_{j=1}^i \bar{p}_{m,j} \right| \left( x_{i+1} - x_i \right) \right) \leq \gamma_m^{1-\mathcal{W}} \right\},
    \label{eq:uc_1w_1d}
\end{equation}
where $x_i < x_j$, $\forall i<j$. The detailed derivation is provided in Online Supplement A.4.

Thirdly, the ambiguity set of the same parametric distribution family is defined as
\begin{equation*}
    \mathcal{P}_m^{\text{Param}} = \left\{\boldsymbol{p}_m \middle\vert \ p_{m,i} = \mathbb{P}\left( X_m = x_i\right), \ \forall i = 1, \dots,  |\Omega|, \ \text{where } X_m \sim \mathcal{D}_m(\theta_m), \ \forall \theta_m \in \Theta_m \right\},
\end{equation*}
where $X_m$ denotes the R.V. for the $m$th input model, $\mathcal{D}(\theta_m)$ is a certain member within a pre-specified distribution family with its parameter $\theta_m$, and $\Theta_m$ is the set of candidate parameters. Here, the magnitude $|\Theta_m|$ of the range in which the parameter varies can be interpreted as the ambiguity set size parameter $\gamma_m$. 

For example, if $X_m$ is a binomial R.V. with parameters $(N_m^{\text{Bin}}$, $p_m^{\text{Bin}}$), the ambiguity set can be expressed as follows:
\begin{equation}
    \mathcal{P}_m^{\text{Bin}} = \left\{\boldsymbol{p}_m \middle\vert \ p_{m,i} = \binom{N_m^{\text{Bin}}}{x_i} \left(p_m^{\text{Bin}}\right)^{x_i} \left(1-p_m^{\text{Bin}}\right)^{N_m^{\text{Bin}}-x_i}, \ \forall i = 1,\dots, |\Omega|, \  \forall (N_m^{\text{Bin}}, p_m^{\text{Bin}}) \in \Theta_m \right\} .
    \label{eq:uc_binomial}
\end{equation}

As another example, let us consider a discretized version of the Rayleigh distribution family. This set will be used in our case study that analyzes wind turbine simulator outputs. By letting the input R.V.'s probability mass be proportional to the probability density of Rayleigh distribution with an input shift, we get the following ambiguity set. 
\begin{equation}
    \mathcal{P}_m^{\text{Rayleigh}} = \left\{\boldsymbol{p}_m \middle\vert \ p_{m,i} \propto \frac{x_i - \Delta_m}{\left(\sigma_m^{\text{Rayleigh}}\right)^2} e^{-\frac{1}{2}\left( \frac{x_i - \Delta_m}{\sigma_m^{\text{Rayleigh}}} \right)^2}, \ \forall i = 1, \dots, |\Omega|, \ \forall (\sigma_m^{\text{Rayleigh}}, \Delta_m) \in \Theta_m \right\},
    \label{eq:uc_rayleigh_with_shift}
\end{equation}
with the pair of Rayleigh scale parameter $\sigma_m^{\text{Rayleigh}}$ and input shift $\Delta_m$. We note that the choice of the parametric family is not limited to the examples here; one may instead select any other family depending on prior domain expertise.

Finally, we define the ambiguity set based on distribution moments as follows:
\begin{equation}
    \begin{aligned}
        &\mathcal{P}_m^{\text{Moment}} \\
        &= \left\{ \boldsymbol{p}_{m} \middle\vert 
        \begin{array}{l} 
            \left( \sum\limits_{i=1}^{|\Omega|} p_{m,i} x_i - \bar{\boldsymbol{\mu}}_m\right)^{\text{T}} \bar{\boldsymbol{\Sigma}}_m^{-1} \left( \sum\limits_{i=1}^{|\Omega|} p_{m,i} x_i - \bar{\boldsymbol{\mu}}_m \right) \leq \gamma_{1,m} \\
            \sum\limits_{i=1}^{|\Omega|} p_{m,i} \left( x_i -\bar{\boldsymbol{\mu}}_m\right) \left( x_i -\bar{\boldsymbol{\mu}}_m\right)^{\text{T}} \preceq \gamma_{2,m}^{ub} \bar{\boldsymbol{\Sigma}}_m \\ 
            \sum\limits_{i=1}^{|\Omega|} p_{m,i} \left( x_i -\bar{\boldsymbol{\mu}}_m\right) \left( x_i -\bar{\boldsymbol{\mu}}_m\right)^{\text{T}}\succeq 
            \gamma_{2,m}^{lb} \bar{\boldsymbol{\Sigma}}_m 
            + 2 \left( \sum\limits_{i=1}^{|\Omega|} p_{m,i} x_i - \bar{\boldsymbol{\mu}}_m\right)\left( \sum\limits_{i=1}^{|\Omega|} p_{m,i} x_i - \bar{\boldsymbol{\mu}}_m\right)^T
        \end{array}\right\},
    \end{aligned}
    \label{eq:uc_moment_pmf}  
\end{equation}
where $\bar{\boldsymbol{\mu}}_m$ and $\bar{\boldsymbol{\Sigma}}_m$ are the mean vector and covariance matrix of the nominal input vector $\bar{X}_m$, respectively, and $\gamma_{1,m}, \gamma_{2,m}^{lb},$ and $ \gamma_{2,m}^{ub}$ are positive scalar values determining the level of uncertainty. This set is an extension of the ambiguity set proposed in~\cite{delage2010distributionally}. The original set in~\cite{delage2010distributionally} bounds above the first and second-order moments, but we also include the third constraint to further limit the second moment to be bounded below. This lower bound is included because, in the problem under consideration in this study, both extreme instances can result in the largest estimator variance. Online Supplement A.5 discusses how we construct this new ambiguity set in detail.

With $\gamma_{1,m} = 0$ and $\gamma_{2,m}^{lb} = \gamma_{2,m}^{ub} = 1$, this ambiguity set consists of the distributions which have the same first and second moments as the nominal distribution. But, this does not imply that the ambiguity set includes the nominal distribution only. 

\subsection{Solving DR-Strat Problem}
\label{subsec:solving_DR-Strat_problem}

This section discusses how to solve the DR-strat problem in \eqref{eq:min_max_opt_model}. In our case, the variable to be optimized is the sampling vector $\boldsymbol{n}$, which is the decision vector in the outer problem with regard to the inner maximization problem's objective value $v(\boldsymbol{n})$. For calculating  $v(\boldsymbol{n})$ given the sampling vector $\boldsymbol{n}$, one can either apply the iterative algorithm or use a nonlinear solver. In this study, we utilize open-source solvers with implementation details provided in Online Supplement B.

The challenge lies in solving the outer problem. One may consider enumerating all potential candidates $\boldsymbol{n}$'s and choosing the one that generates the smallest  $v(\boldsymbol{n})$. This naive approach is, however, not computationally efficient, even if it is possible to compute. The number of possible solutions, ${}_{N_T-1}C_{N_T-K}$ by the formula of combination with repetition, becomes extremely huge (e.g., approximately $10^{10}$ for $N_T=100$ and $K=7$) even with moderate $N_T$ and $K$, making an exhaustive search computationally intractable. 

In the literature, several algorithms have been presented for solving a bi-level optimization problem (e.g., using a single-level reduction or KKT conditions)~\citep{bard2013practical}. Unfortunately, the objective function of our inner problem is a non-convex form of the pmf $\boldsymbol{p}$, preventing us from employing existing techniques. Recent studies on robust decision-making show that evolutionary approaches, such as a genetic algorithm, can be used to solve analytically intractable problems, but they tend to heavily focus on exploitation. 

We utilize BO, a probabilistic global optimization approach which is known to strike a balance between exploration and exploitation,  and to be effective in handling multi-local-optima~\citep{snoek2012practical}. BO models the variable-function value relationship with GP, which iteratively updates as new observations become available. Specifically, we start with an initial set $\mathcal{D}_{sv}$ of sampling vectors and the corresponding set $\mathcal{V}_{inner} = \{v(\boldsymbol{n}), \forall \boldsymbol{n} \in \mathcal{D}_{sv}\}$ of the objective values of the inner problem. Then, we model the relationship between the sampling vector and its corresponding objective value with GP. A new candidate sampling vector  $\boldsymbol{n}^{new}$ is determined by maximizing the acquisition function (ACQ). Among various ACQs, we utilize the following expected improvement over the best objective value found so far.
\begin{equation}
    \text{EI}(\boldsymbol{n}) = \mathbb{E} \left[ \max \left(v(\boldsymbol{n}^{best}) - v(\boldsymbol{n}), 0 \right) \right], 
    \label{eq:EI}
\end{equation}
which can be calculated using the mean and 
variance of the GP posterior at $\boldsymbol{n}$. 

This new sampling vector $\boldsymbol{n}^{new}$ and its objective value $v(\boldsymbol{n}^{new})$ are added to $\mathcal{D}_{sv}$ and $\mathcal{V}_{inner}$, respectively. If the new objective value is better (lower) than the current best, $\boldsymbol{n}^{new}$ replaces $\boldsymbol{n}^{best}$. These steps are repeated until a stopping criterion is met. During the iteration, we allow the elements of $\boldsymbol{n}$ to have continuous values rather than restricting them to integers. When the iteration finally terminates, we round up the obtained $\boldsymbol{n}^{best}$. Additional details are provided in Online Supplement B.

\section{Numerical Experiments}
\label{sec:numerical_experiments}

This section assesses the effectiveness of the DR-strat method. Section~\ref{subsec:benchmark_model} describes a modified stratified sampling approach as the benchmark model that takes into account multiple input models without uncertainties. We implement the proposed methodology and compare it with the benchmark model in two experimental settings: a numerical example in Section~\ref{subsec:toy_example} and the case study involving wind turbine reliability in Section~\ref{subsec:case_study}. 
Online Supplement D.4 also provides additional experimental results with two-dimensional input.

\subsection{Benchmark Model}
\label{subsec:benchmark_model}

Section~\ref{subsec:recap_stratification} has outlined the conventional stratified sampling method, which handles a single input model. To the best of our knowledge, no prior studies in stratified sampling consider multiple distributions with input uncertainty. For fair comparison, we use a modified approach as our benchmark model that ignores input uncertainty while handling multiple input models. Specifically, assuming complete information about input models, the benchmark model treats the nominal distributions as true input models. Similar to the proposed DR-strat, it uses a single reference distribution during the sampling phase and then estimates the response for each input model using~\eqref{eq:dr_strat_estimator}. With the goal of obtaining a sampling vector that minimizes the maximum estimator variance among multiple nominal input models, it formulates the following problem. 
\begin{equation*}
    \begin{aligned}
        \text{(Str-M)}\ \underset{\boldsymbol{n}}{\min} \ \underset{1 \leq m \leq M}{\max} \ & \sum_{k=1}^{K} \frac{1}{n_k} \left( \left( \omega_{ref,k} \sum_{ i \in I_k } \mathbb{E}\left[ g(x_i) \right] \frac{ \bar{p}_{m,i}^2}{ p_{ref,i} } \right) - \left( \sum_{ i \in I_k} \mathbb{E}\left[ g(x_i) \right] \bar{p}_{m,i} \right)^2 \right)\\
        s.t. \ &\sum_{k=1}^K n_k = N_T.
    \end{aligned}
\end{equation*}
Please note that the objective term in the optimization problem (Str-M) does not have a maximum operator $\underset{\boldsymbol{p}_m \in \mathcal{P}_m}{\max}$ which reflects the uncertainty in the $m$th input model, unlike that in (DR-Str) in \eqref{eq:min_max_opt_model}. Let $\boldsymbol{n}^{\text{Str-M}}$ denote the optimal sampling vector of (Str-M), where  $\text{M}$ in the superscript implies the consideration of \textit{multiple} input models, in contrast to $\boldsymbol{n}^{\text{Str}}$ in Section~\ref{subsec:recap_stratification}. 

\subsection{Toy Example}
\label{subsec:toy_example}

\subsubsection{Experimental Setting.} 
\label{subsubsec:toy_example_experiment_setting}

Consider estimating the tail probability $\mathbb{P}(Y(X)>l)$ with a one-dimensional input $X$. Mimicking the standard normal input R.V. in the example in~\cite{myoung2022optimal}, we employ the following scaled binomial R.V.s $\bar{X}_1$ and $\bar{X}_2$ as the nominal distributions of two input models, with the domain of $\bar{B}_1$ and $\bar{B}_2$ as $\{23, 24, \dots, 57\}$.
\begin{equation*}
    \bar{X}_1 = \frac{\bar{B}_1- (80\times0.5)}{\sqrt{80\times 0.5^2}}, \ \bar{X_2} = \frac{\bar{B_2}- (80\times0.5)}{\sqrt{80\times 0.5^2}}, \ \text{ where } \bar{B}_1 \sim \text{Bin}(75, 0.55), \ \bar{B}_2 \sim \text{Bin}(85, 0.45). 
\end{equation*}

For the output model, we use the same model in~\cite{myoung2022optimal} and define the conditional output given a certain input to be $Y|\{X=x\} \sim \mathcal{N}\left(\mu_Y(x), \sigma_Y(x)\right)$ with
\begin{equation}
    \begin{gathered}
    \mu_Y(x) = 0.95x^2 (1 + 0.5\cos(10x) + 0.5\cos(20x)),\\
    \sigma_Y(x) = 1 + 0.7|x| + 0.4\cos(x) + 0.3\cos(14x).
    \end{gathered}
    \label{eq:toy_output}
\end{equation}
To meet the target performance measure (the tail probability) values with the nominal distributions as $\mathbb{P}(Y(\bar{X}_1)>l) = 0.0428$ and $\mathbb{P}(Y(\bar{X}_2)>l) = 0.0564$, a threshold $l$ is set to be $5.2$.

For the input domain $\Omega = \left\{ x_i | x_i = (i-40)/\sqrt{20}, \ \forall i = 23, \dots, 57 \right\}$, we consider 7 strata $S_k = \left\{x_i | x_i = (i-40)/\sqrt{20}, \ \forall i = 23+5(k-1), \dots, 22+5k \right\}$ for $k = 1, \dots, 7$. The total simulation budget $N_T$ is 100. We use the mean of the two nominal distributions as the reference distribution for initial sampling ( i.e., $\mathbb{P}\left( X_{ref} = x_i \right) = \left( \mathbb{P}\left( \bar{X}_1 = x_i \right) + \mathbb{P} \left( \bar{X}_2 = x_i \right) \right)/2, \ \forall x_i \in \Omega$ ). We recommend choosing the reference distribution near the nominal distributions which the ambiguity sets are constructed around. Obtaining the optimal reference distribution remains a subject of our future study. Further, we assume that conditional output means $\{\mathbb{E}[g(x_i)]\}_{i =1}^{|\Omega|}$ are known as in~\eqref{eq:toy_output}. In reality, we can estimate them via learning a meta-model with the results obtained from running the pilot stage simulations. To construct the four types of ambiguity sets, we utilize \eqref{eq:uc_l2_1d}, \eqref{eq:uc_1w_1d}, \eqref{eq:uc_binomial}, and \eqref{eq:uc_moment_pmf}. The detailed settings for set size parameters are provided in Online Supplement C.1.

\subsubsection{Instances within Ambiguity Sets.}    \label{subsubsec:toy_example_instances_within_uncertainty_sets}

\begin{figure}[t!]
    \centering  
    \begin{subfigure}[b]{0.26\textwidth}
        \centering
        \includegraphics[width=\textwidth]{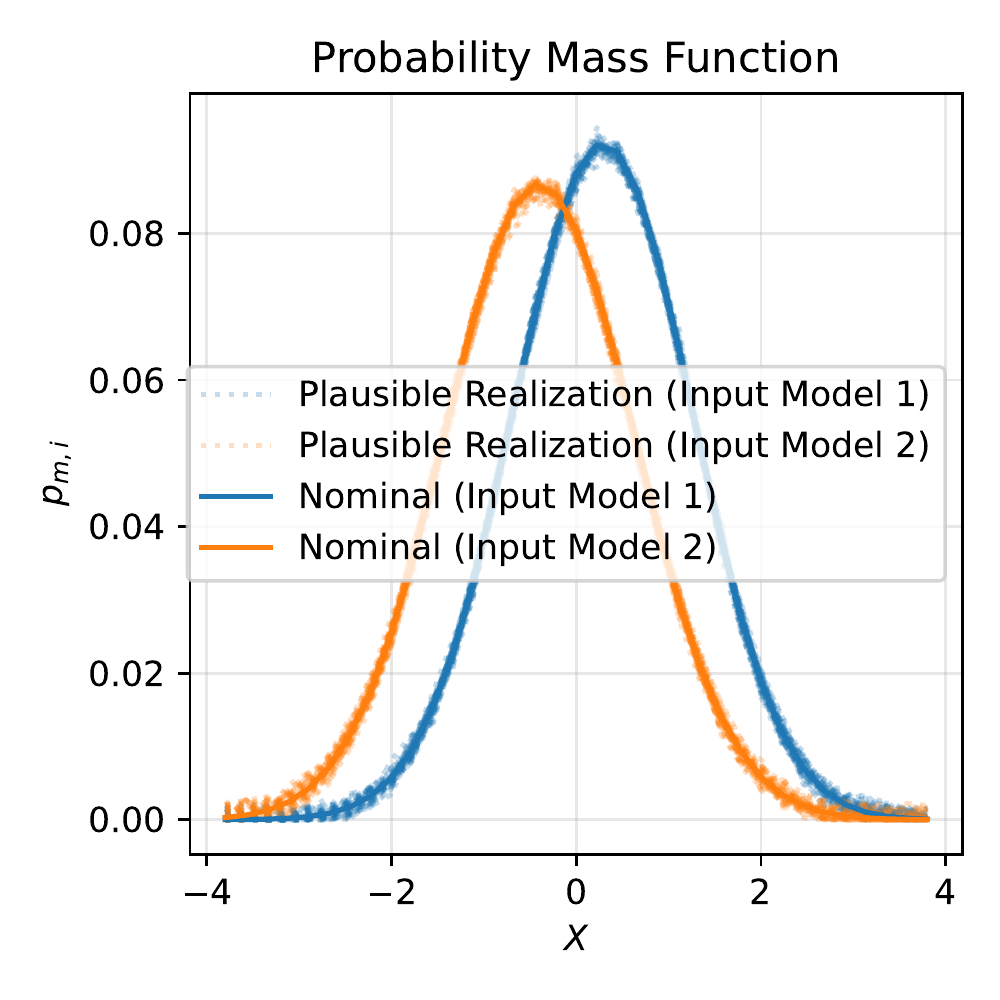}
        \caption{$L_2$-Norm}
        \label{subfig:instances_within_uc_l2_toy}
    \end{subfigure}
    \begin{subfigure}[b]{0.26\textwidth}
        \centering
        \includegraphics[width=\textwidth]{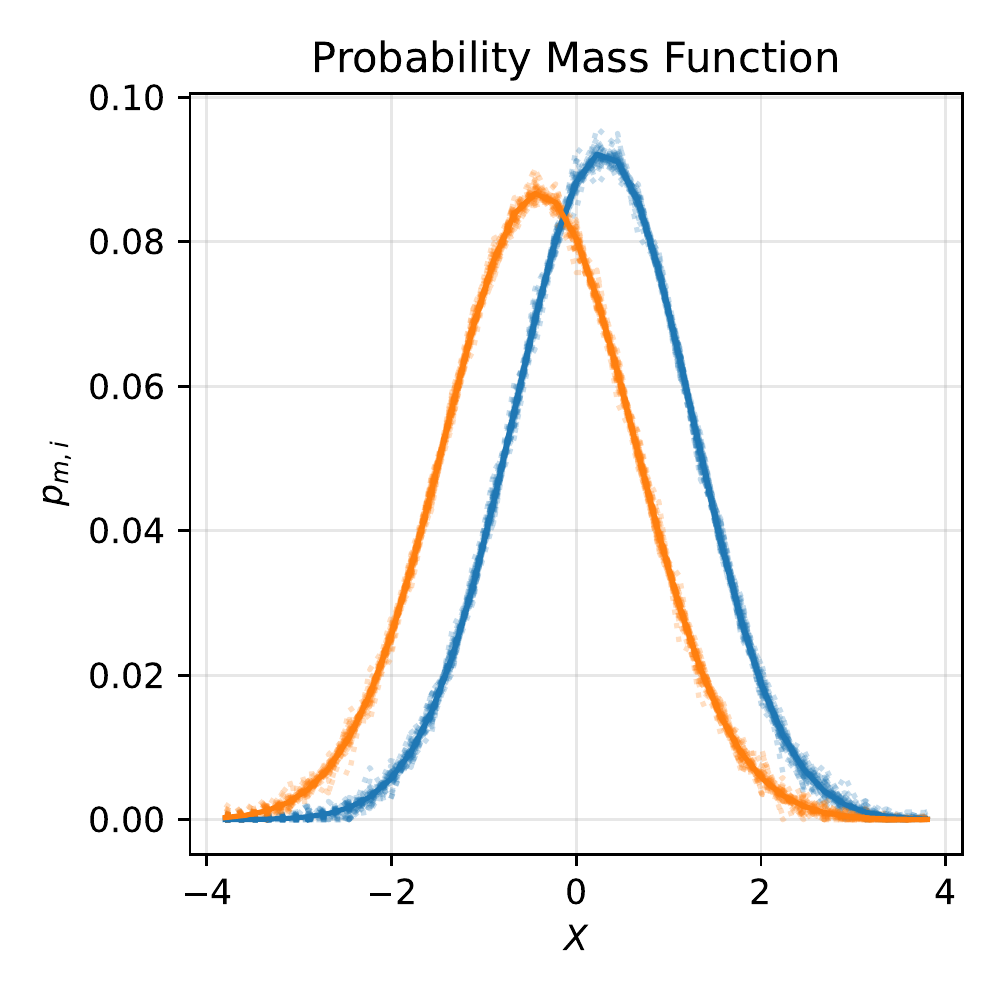}
        \caption{$1$-Wasserstein Distance}
        \label{subfig:instances_within_uc_1W_toy}
    \end{subfigure}
    \\
    \begin{subfigure}[b]{0.26\textwidth}
        \centering
        \includegraphics[width=\textwidth]{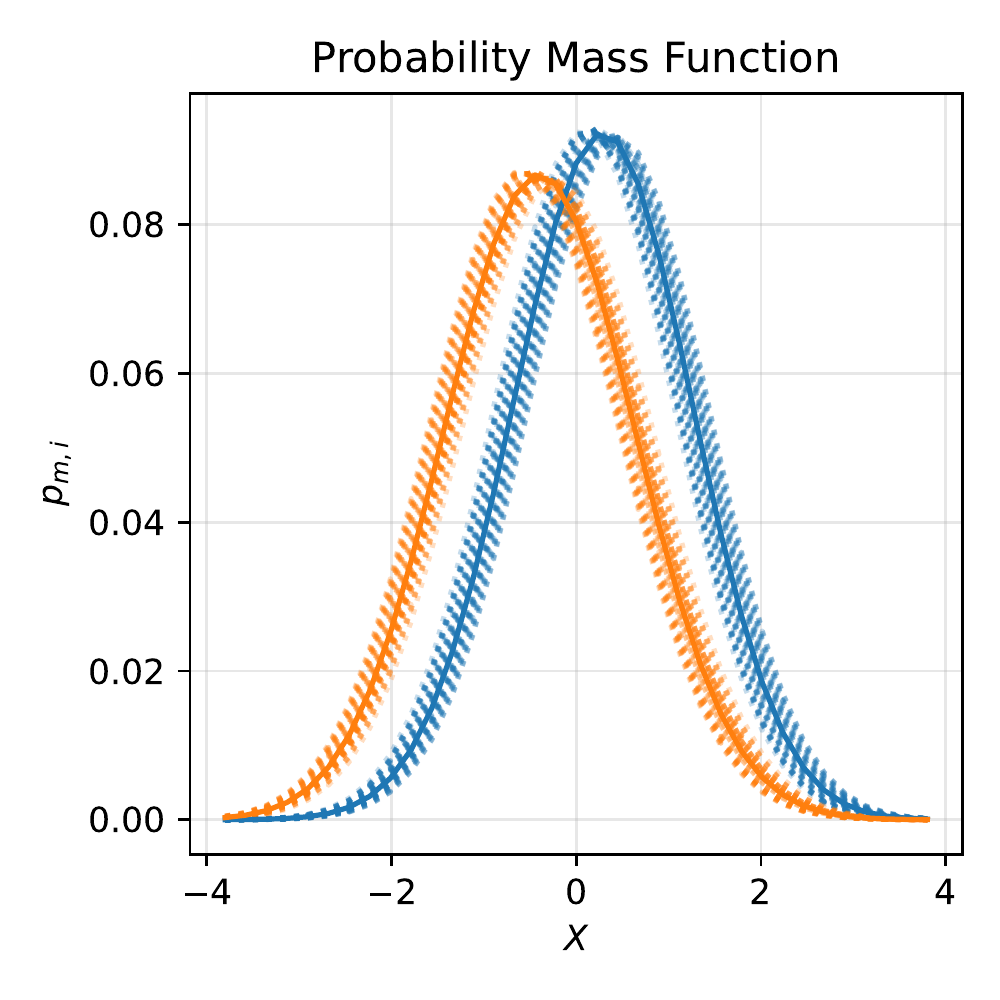}
        \caption{Parametric Family}
        \label{subfig:instances_within_uc_param_toy}
    \end{subfigure}
    \begin{subfigure}[b]{0.26\textwidth}
        \centering
        \includegraphics[width=\textwidth]{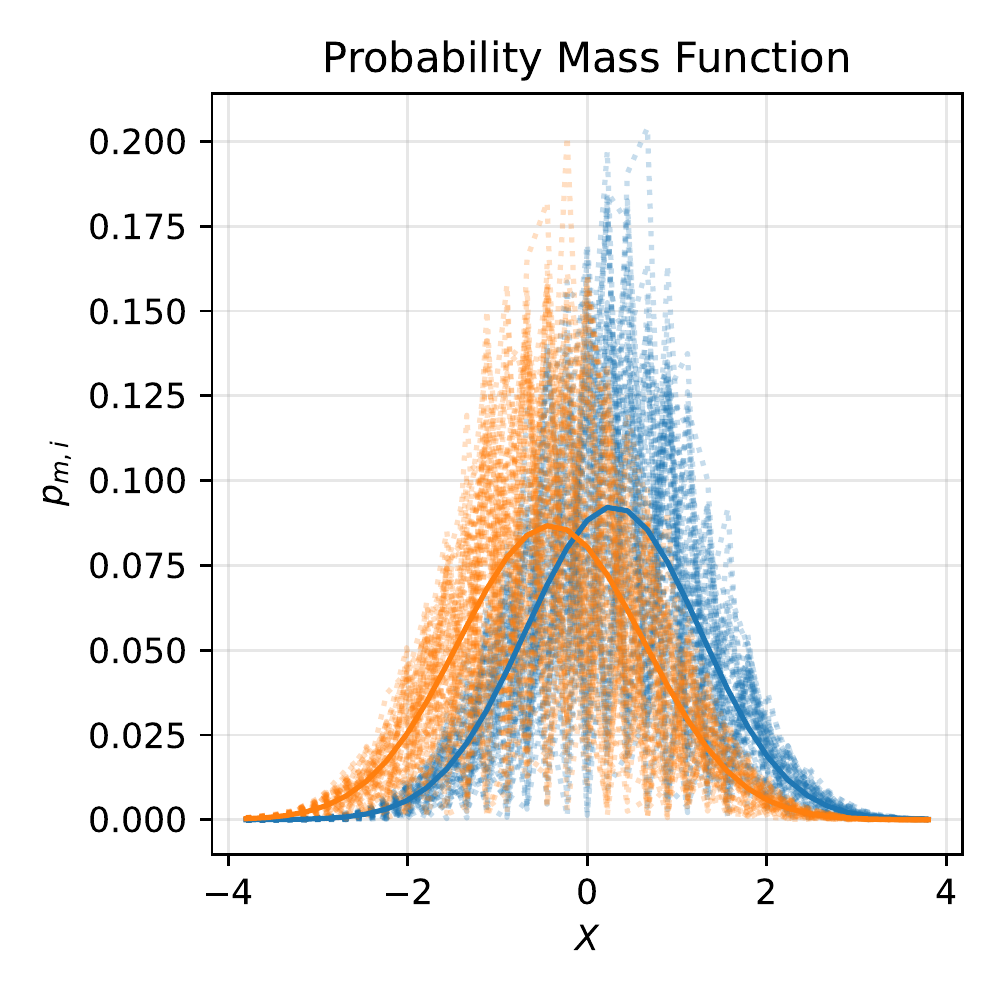}
        \caption{Distribution Moments}            \label{subfig:instances_within_uc_moments_toy}
    \end{subfigure}
    \caption{Plausible Input Distributions within Four Ambiguity Sets (Note: Although the input distributions are discrete, we exhibit continuous patterns for demonstration purposes.)}
    \label{fig:plausible_input_distributions_within_UCs_toy}
\end{figure}

We first depict plausible distributions in each ambiguity set in Figure~\ref{fig:plausible_input_distributions_within_UCs_toy} for the two input models. Solid and dotted curves represent the pmfs of the nominal and plausible distributions, respectively. Depending on the underlying similarity measure, the plausible distributions show different shapes. 

The pmfs in the ambiguity sets constructed with the two discrepancy measures, shown in Figures~\ref{subfig:instances_within_uc_l2_toy} and \ref{subfig:instances_within_uc_1W_toy}, similarly show moderate spikes. They are alike to many empirical distributions fitted from historical  
data. Still, the differences between the two sets exist. In the $1$-$\mathcal{W}$ distance-based set, plausible distributions tend to differ from the nominal distribution at fewer input points, compared to those in $L_2$-norm-based set. This is because $L_2$-norm calculates the difference of pmfs individually at each input value in~\eqref{eq:uc_l2_1d}, whereas the $1$-$\mathcal{W}$ distance is calculated in a cumulative manner in~\eqref{eq:uc_1w_1d}.

Next, Figure~\ref{subfig:instances_within_uc_param_toy} shows smooth pmfs in the ambiguity set of a parametric family. 
This set type is desirable if historical data is to be fitted to a pre-specified parametric distribution. Finally, pmfs in the moment-based ambiguity set in Figure~\ref{subfig:instances_within_uc_moments_toy} depict the most jagged shapes. These spiky pmfs appear because the set constraints restrict only the first two moments but not the distribution shape. As a given input's probability mass can vary greatly, feasible realizations may deviate dramatically from the nominal distribution pattern. Employing this set type may result in excessive conservatism when incorporating the input model uncertainty. As a result, the moment-based set should be used only when the true distribution possibly has an unusual pmf form.

\subsubsection{Implementation Results.}
\label{subsubsec:sampling_vector_comparison}

\begin{figure}[t!]
    \centering  
    \begin{subfigure}[b]{0.26\textwidth}
        \centering
        \includegraphics[width=\textwidth]{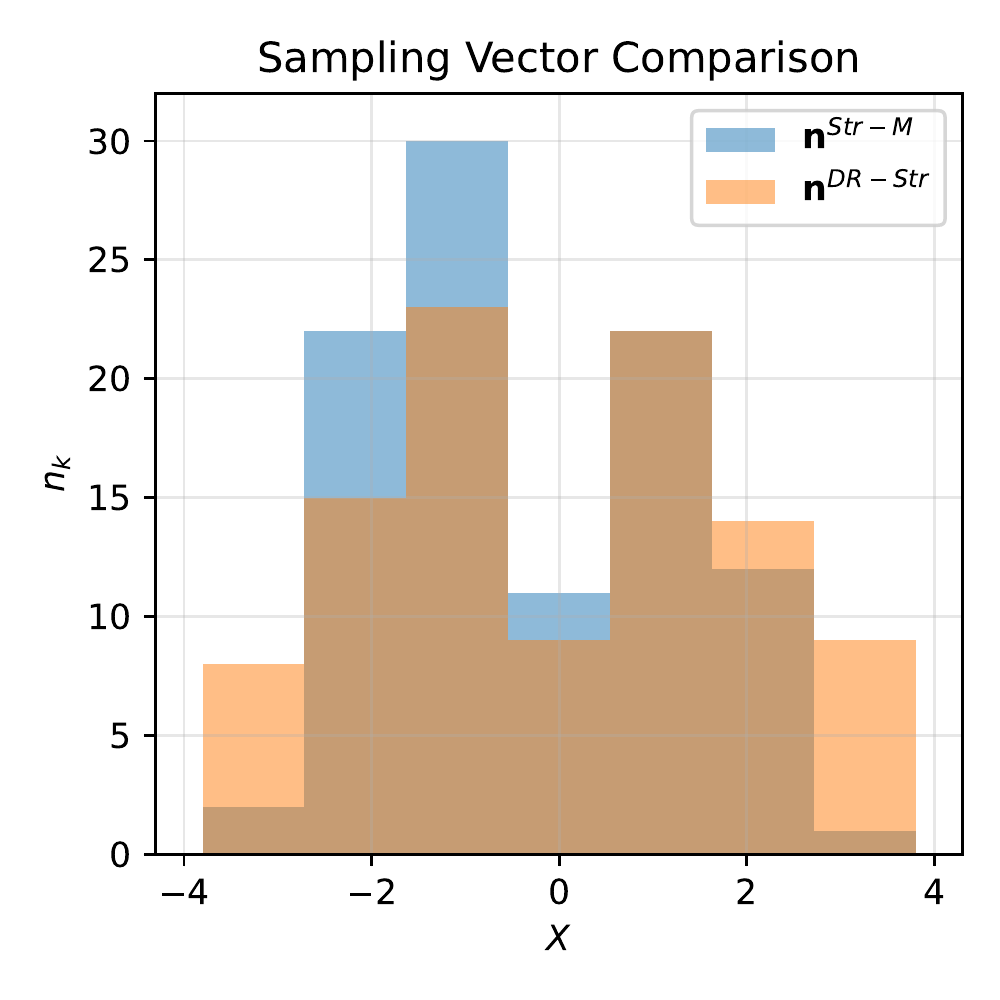}
        \caption{$L_2$-Norm}
        \label{subfig:stratification_vector_comparison_l2_toy}
    \end{subfigure}
    \begin{subfigure}[b]{0.26\textwidth}
        \centering
        \includegraphics[width=\textwidth]{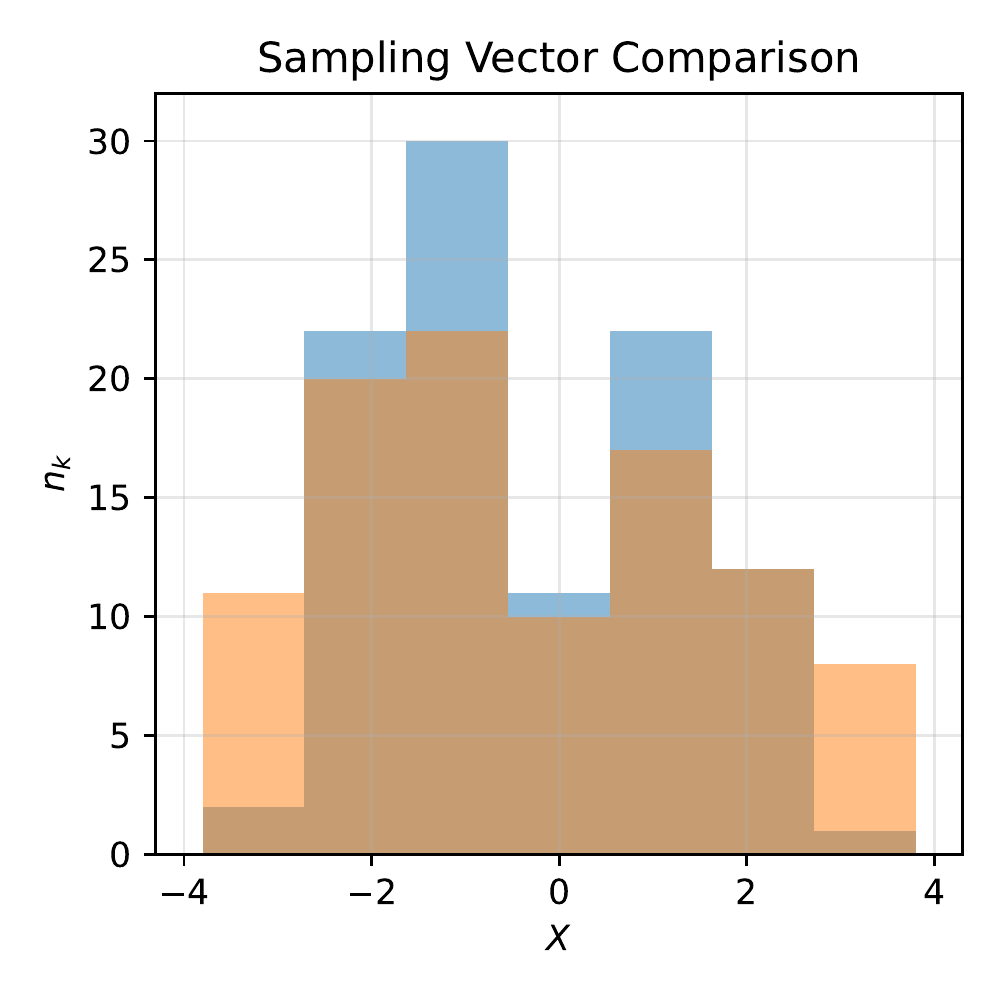}
        \caption{$1$-Wasserstein Distance}
        \label{subfig:stratification_vector_comparison_1W_toy}
    \end{subfigure}
    \\
    \begin{subfigure}[b]{0.26\textwidth}
        \centering
        \includegraphics[width=\textwidth]{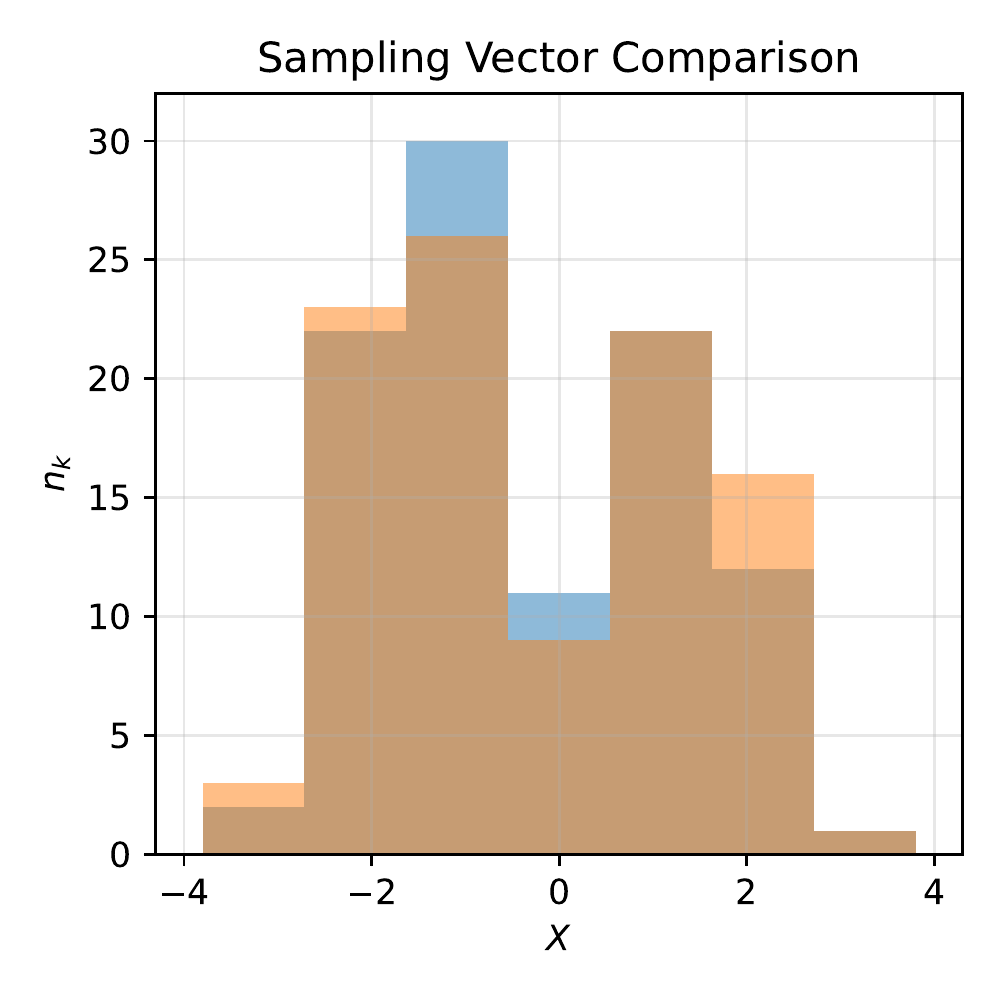}
        \caption{Parametric Family}
        \label{subfig:stratification_vector_comparison_param_toy}
    \end{subfigure}
    \begin{subfigure}[b]{0.26\textwidth}
        \centering
        \includegraphics[width=\textwidth]{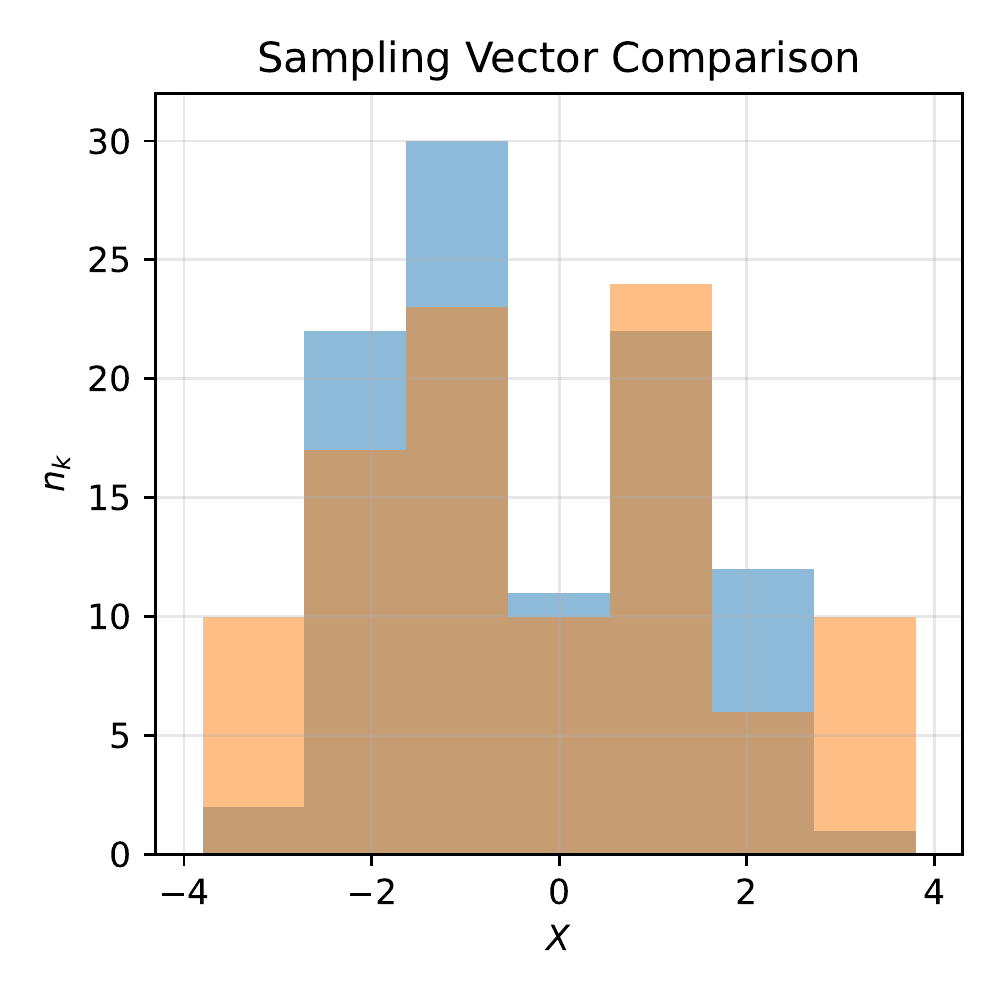}
        \caption{Distribution Moments}            \label{subfig:stratification_vector_comparison_moments_toy}
    \end{subfigure}
    \caption{Comparison of Sampling Vectors Obtained from DR-Strat and Benchmark Model}            \label{fig:stratification_vector_comparison_Strat_vs_DR_Strat_toy}
\end{figure}

We compare the sampling vectors from DR-strat and the benchmark method, using the same reference (sampling) distribution $F_{ref}$ for both methods. Figure~\ref{fig:stratification_vector_comparison_Strat_vs_DR_Strat_toy} depicts $\boldsymbol{n}^{\text{DR-Str}}$ and $\boldsymbol{n}^{\text{Str-M}}$ for each ambiguity set design. We observe that $\boldsymbol{n}^{\text{Str-M}}$ focuses intensively on particular strata (allocating $74\%$ of the total simulation budget the $2$nd, $3$rd, and $5$th strata) where both the input probability and conditional output variance $\text{Var} \left[g(X)|X \in S_k \right]$ are relatively high. On the other hand, $\boldsymbol{n}^{\text{DR-Str}}$ tends to be more stretched even to strata with low probability but high conditional output variance. In all set designs, $\boldsymbol{n}^{\text{DR-Str}}$ has a smaller maximum ($\underset{1 \leq k \leq K}{\max}  n_k^{\text{DR-Str}}$) and a larger minimum budget ($\underset{1 \leq k \leq K}{\min}  n_k^{\text{DR-Str}}$), compared to $\boldsymbol{n}^{\text{Str-M}}$, indicating the conservative tendency of the proposed method.

\begin{figure}[t!]
    \centering  
    \begin{subfigure}[b]{0.26\textwidth}
        \centering
        \includegraphics[width=\textwidth]{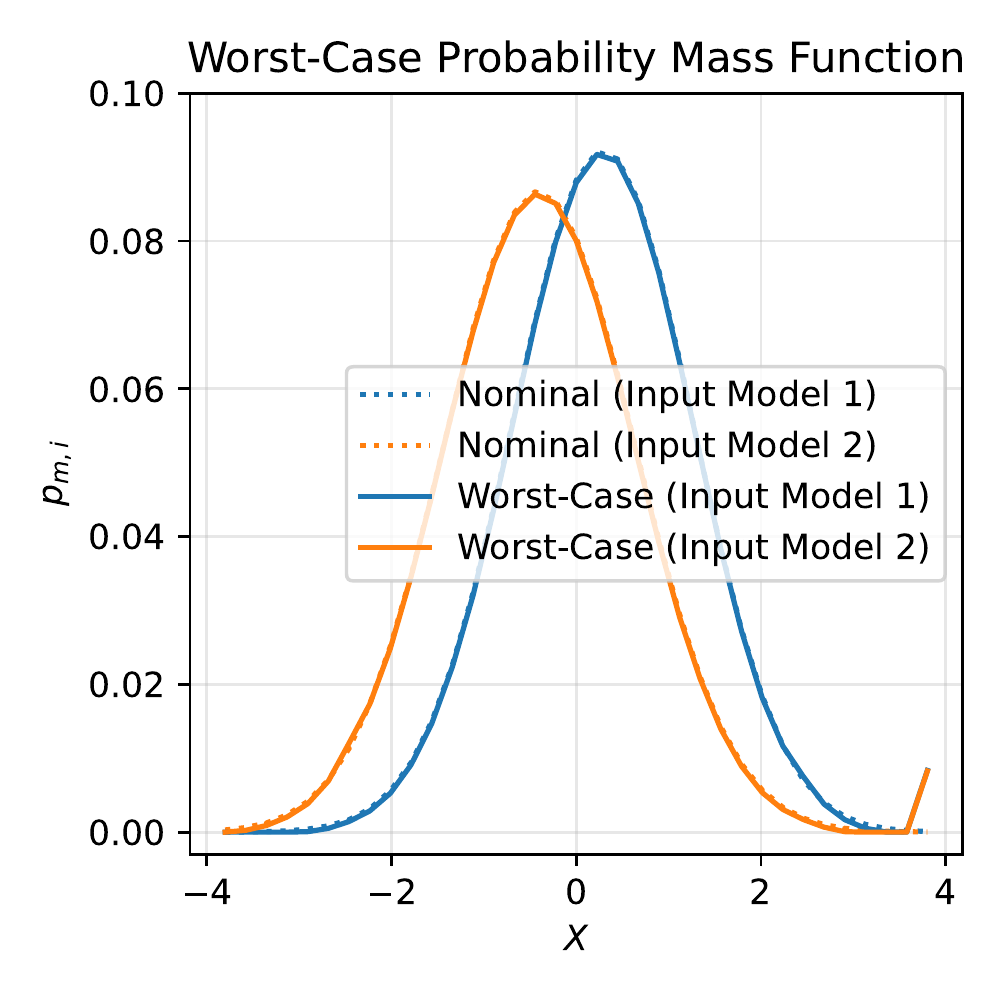}
        \caption{$L_2$-Norm}
        \label{subfig:worst_case_distributions_within_UCs_l2_toy}
    \end{subfigure}
    \begin{subfigure}[b]{0.26\textwidth}
        \centering
        \includegraphics[width=\textwidth]{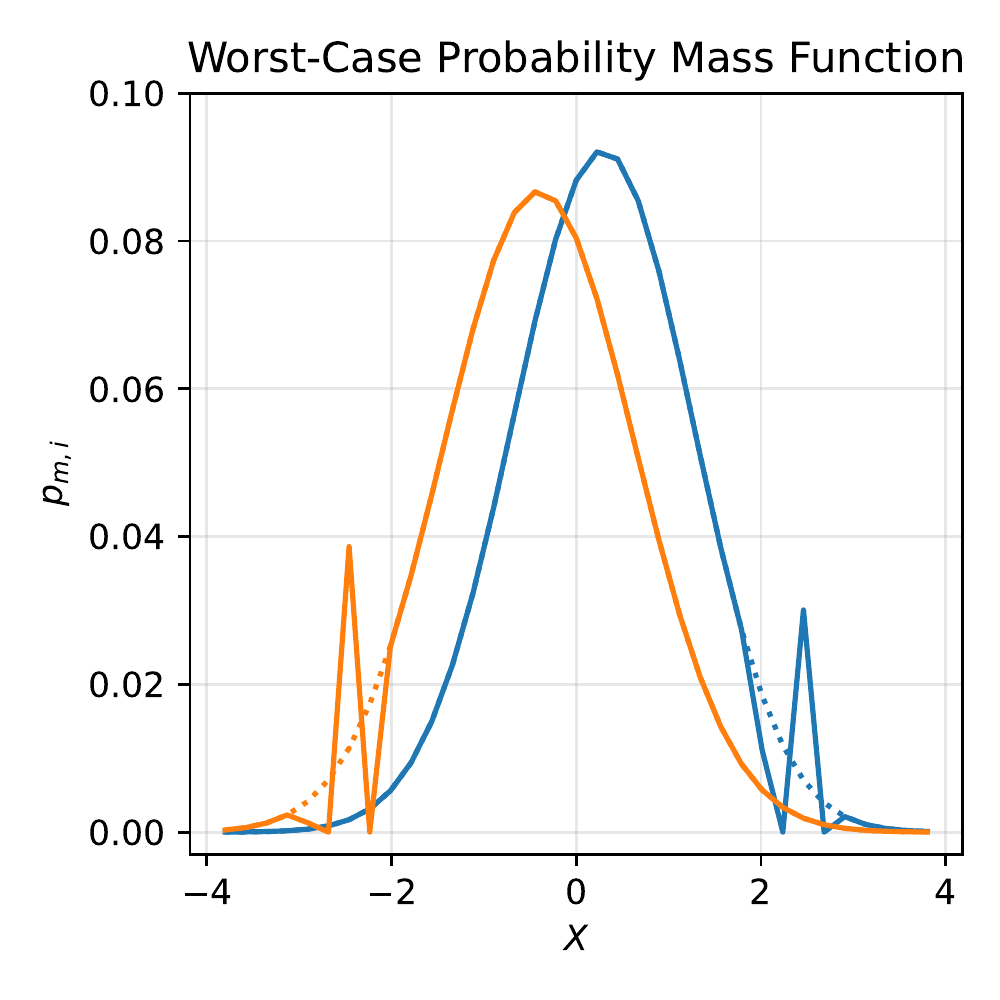}
        \caption{$1$-Wasserstein Distance}
        \label{subfig:worst_case_distributions_within_UCs_1W_toy}
    \end{subfigure}
    \\
    \begin{subfigure}[b]{0.26\textwidth}
        \centering
        \includegraphics[width=\textwidth]{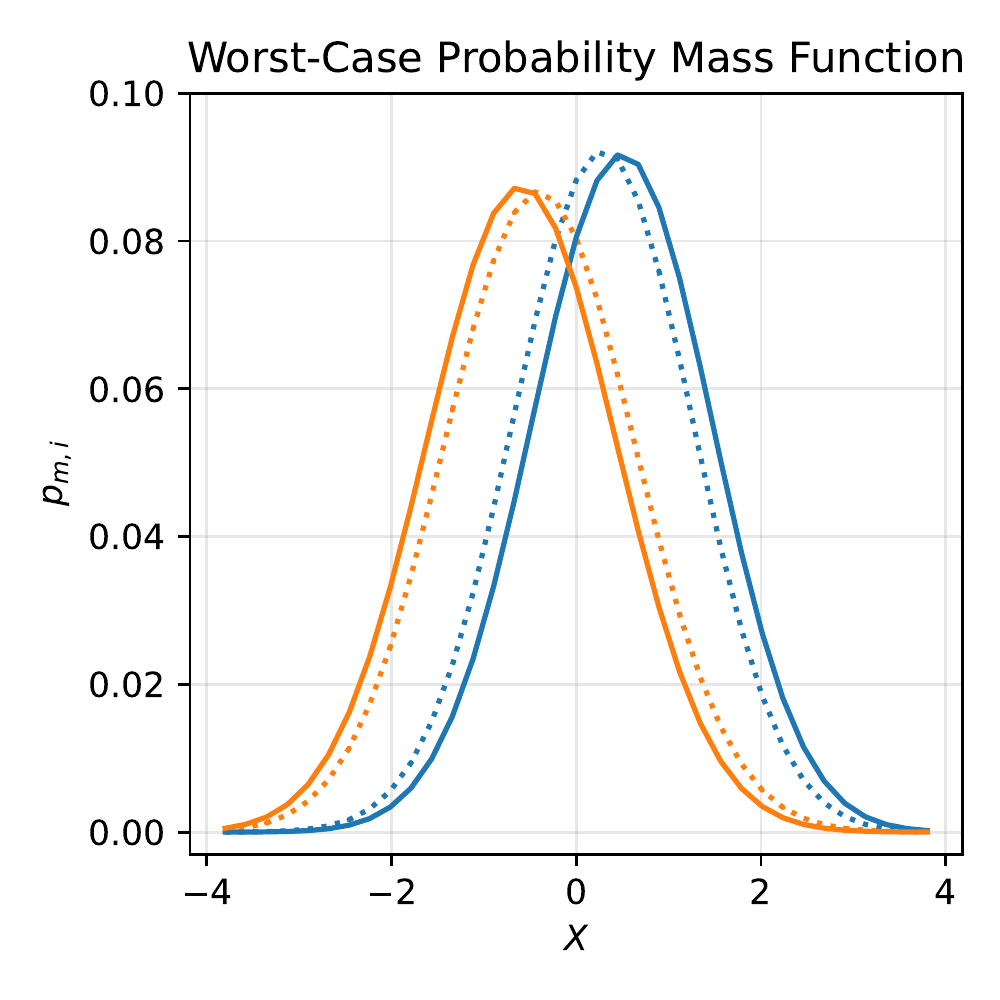}
        \caption{Parametric Family}
        \label{subfig:worst_case_distributions_within_UCs_param_toy}
    \end{subfigure}
    \begin{subfigure}[b]{0.26\textwidth}
        \centering
        \includegraphics[width=\textwidth]{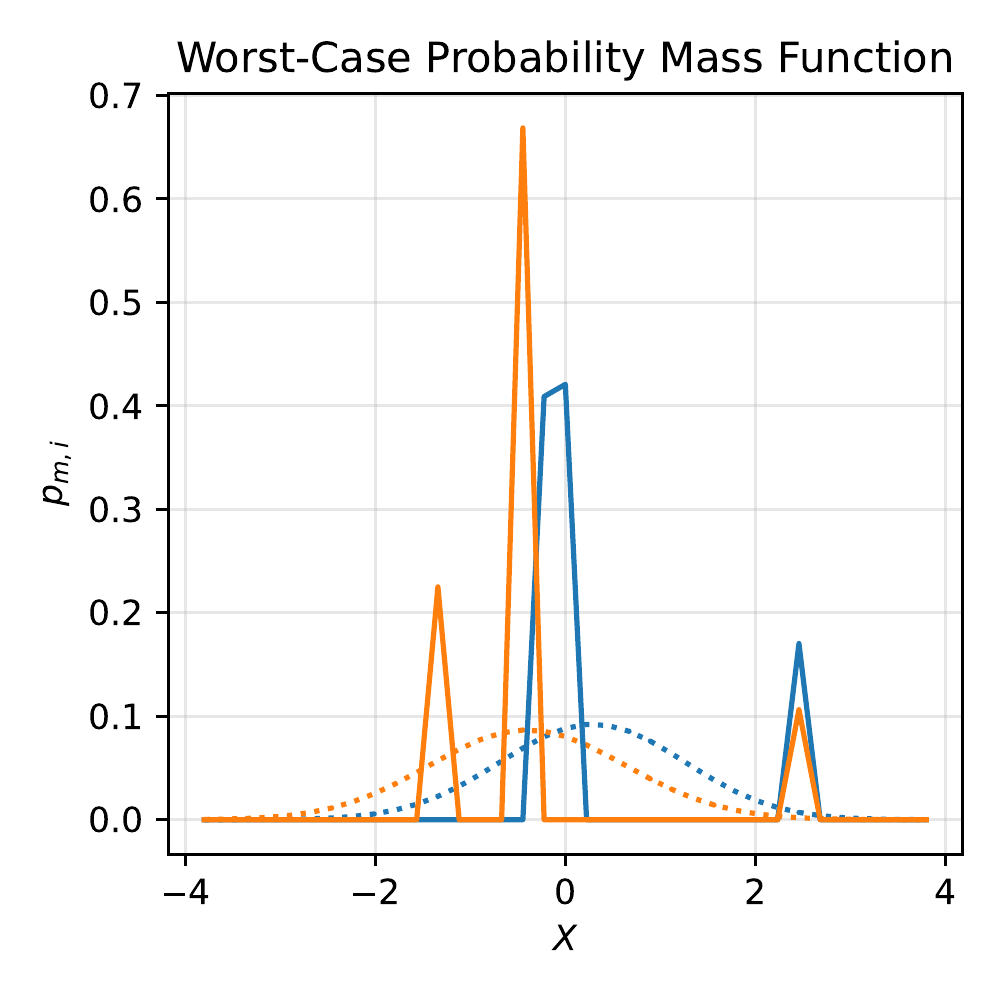}
        \caption{Distribution Moments}
        \label{subfig:worst_case_distributions_within_UCs_moments_toy}
    \end{subfigure}
    \caption{Worst-Case Distributions within Four Ambiguity Sets}
    \label{fig:worst_case_distributions_within_UCs_toy}
\end{figure}

The resulting sampling vector for each ambiguity set is further investigated in conjunction with the corresponding worst-case distribution in Figure~\ref{fig:worst_case_distributions_within_UCs_toy}.
The vector $\boldsymbol{n}^{\text{DR-Str}}$ from the two discrepancy-based sets exhibits similar patterns in Figures~\ref{subfig:stratification_vector_comparison_l2_toy} and~\ref{subfig:stratification_vector_comparison_1W_toy}. However, for the set with $L_2$-norm, more sampling budgets are allocated in the last stratum compared to the set with $1$-$\mathcal{W}$ distance. This aligns with that the worst-case distribution of $L_2$-norm occurs in the first input model with the spike near $X=4$ as shown in Figure~\ref{subfig:worst_case_distributions_within_UCs_l2_toy}. On the contrary, the worst-case distribution of $1$-$\mathcal{W}$ distance, shown in Figure~\ref{subfig:worst_case_distributions_within_UCs_1W_toy}, occurs in the second input model with the spike near $X=-2.5$. It drives a higher budget to the corresponding stratum as shown in Figure~\ref{subfig:stratification_vector_comparison_1W_toy}, compared to the allocation in $L_2$-norm. 

Next, Figure~\ref{subfig:worst_case_distributions_within_UCs_param_toy} demonstrates that the worst-case distributions of the parametric set are moved to the side where $|X|$ is large in comparison to the nominal one. Consequently, $\boldsymbol{n}^{\text{DR-Str}}$ in Figure~\ref{subfig:stratification_vector_comparison_param_toy} concentrates more on the strata near $|X|=2$ than $\boldsymbol{n}^{\text{Str-M}}$. Lastly, the sampling vector of the moment-based set in Figure~\ref{subfig:stratification_vector_comparison_moments_toy} shows the most irregular form among the four. While the other three sampling vectors are bimodal, this vector has three modes, similar to the worst-case distributions shown in Figure~\ref{subfig:worst_case_distributions_within_UCs_moments_toy}. Still, its sampling vector is rather smooth, while the worst-case distributions show spikes at certain input points. This is because the inner maximization problem of DR-strat collectively accounts for other conceivable distributions that may have spikes at different input points.

\begin{figure}[t!]
    \centering  
    \includegraphics[width=0.78\textwidth]{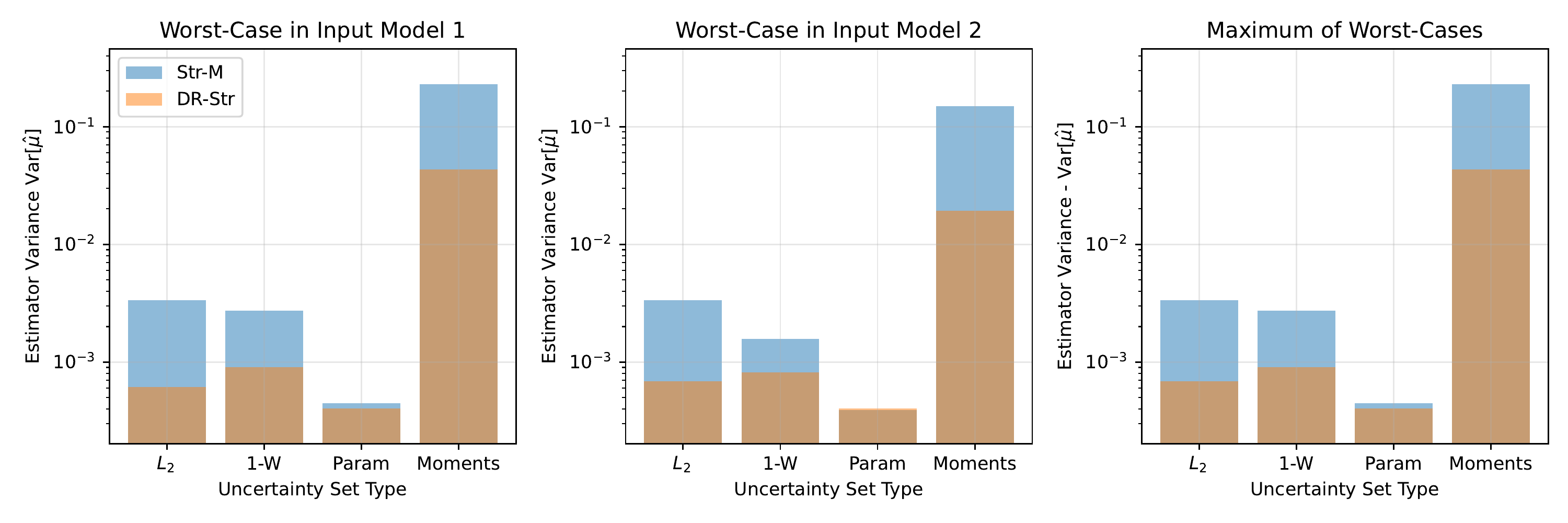}
    \caption{Worst-Case Estimator Variance for Two Input Models and their Maximum Value}
    \label{fig:worst_case_estimator_variances_at_each_timestamp_toy}
\end{figure}

We then compare the estimator variances of DR-strat and the benchmark method, using the derived sampling vectors.
Figure~\ref{fig:worst_case_estimator_variances_at_each_timestamp_toy} depicts the worst-case estimator variance $\underset{\boldsymbol{F}_m \in \mathcal{F}_m}{\max} \text{Var}\left[\hat{\mu}^{\text{DR-Str}}\left(\boldsymbol{n}; F_m\right) \right]$ for $m=1, 2$, as well as their maximum value. The maximum worst-case estimator variance under DR-strat is substantially smaller than the benchmark method's for all four sets, demonstrating its robustness. 

The relative performance of the two methods varies depending on the ambiguity set designs, with the moment-based set exhibiting the most prominent difference, followed by discrepancy-based sets. The magnitude of worst-case variance is considerably larger for the moment-based set than for other sets. Also, DR-strat performs robustly even in circumstances where the true models moderately deviate from the nominal distributions. Online Supplement D.1 showcases such scenarios of the true input model realizations. 

As a final remark, observing the worst-case distributions as well as the realizations with moderate deviations could help determine a proper set design. The moment-based set design tends to consider unrealistically radical distributions and produce overly conservative results. On the other hand, when the true model does not represent the pattern in the same parametric family, a parametric family-based design may produce an overly optimistic set, and the benefit using DR-strat may diminish. The discrepancy-based set designs appear to provide a suitable balance. 

In addition, we conduct sensitivity analysis to assess how the degree of input model uncertainty affects the estimation performance. The true model might deviate from the prediction more (or less) than expected, and the ambiguity set is too small (or large). When compared to the benchmark model, DR-strat leads to lower estimator variance, even when the degree of uncertainty is different from the initial belief, demonstrating its robustness. Online Supplement D.2 provides detailed experimental results and analysis.
 
\subsection{Case Study - Wind Turbine Simulator}
\label{subsec:case_study}

We conduct a case study with a wind turbine simulator. Given a wind condition, the wind turbine simulators—including Turbsim~\citep{jonkman2009turbsim} and FAST~\citep{jonkman2005fast}—generate load responses. Among several load responses, we consider the blade tip defection, which is crucial in analyzing wind turbine reliability~\citep{Choe2016,li2021nonparametric}.  

The simulation input is a 10-min average wind speed. We use the truncated Rayleigh distribution over a support $[3,25]$ with a scale parameter $10\sqrt{2/\pi}$, as recommended in the international standard IEC61400-1~\citep{international2005wind}. We discretize the domain of wind speed into several bins (intervals) in accordance with the widely used binning method in the literature on wind energy. In order to closely mimic the original continuous Rayleigh distribution, we consider a very small bin width of $0.1m/s$. The two input models under consideration have ambiguity sets with the following nominal distributions:
\begin{equation}
 \begin{aligned}
    \bar{p}_{1,i} \propto \frac{x_i - 1.5}{9^2 \times 2/\pi} e^{-\frac{1}{2}\left( \frac{x_i - 1.5}{9\sqrt{2/\pi}} \right)^2}, \ \bar{p}_{2,i} \propto \frac{x_i + 0.5}{11^2 \times 2/\pi} e^{-\frac{1}{2}\left( \frac{x_i + 0.5}{11\sqrt{2/\pi}} \right)^2}, \ \forall i = 1, \dots, |\Omega|,
 \end{aligned}
 \nonumber
\end{equation}
with the domain $\Omega = \{x_i | x_i = 3 + 0.1\times(i-1), \ \forall i = 1, 2, \dots, 220 \}$. We take the average of these two nominal distributions to get the reference distribution. 

\begin{figure}[t!]
    \centering  
    \begin{subfigure}[b]{0.26\textwidth}
        \centering
        \includegraphics[width=\textwidth]{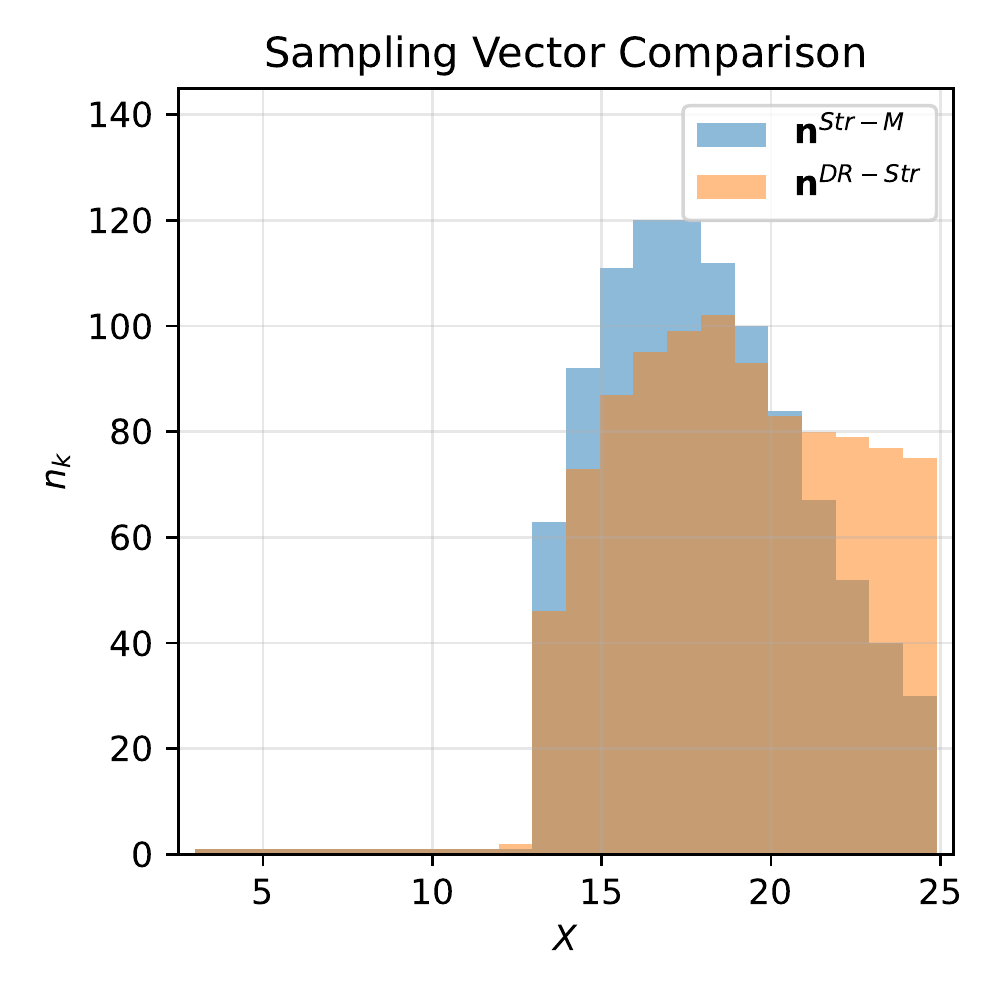}
        \caption{$L_2$-Norm}
        \label{subfig:stratification_vector_comparison_l2_cs}
    \end{subfigure}
    \begin{subfigure}[b]{0.26\textwidth}
        \centering
        \includegraphics[width=\textwidth]{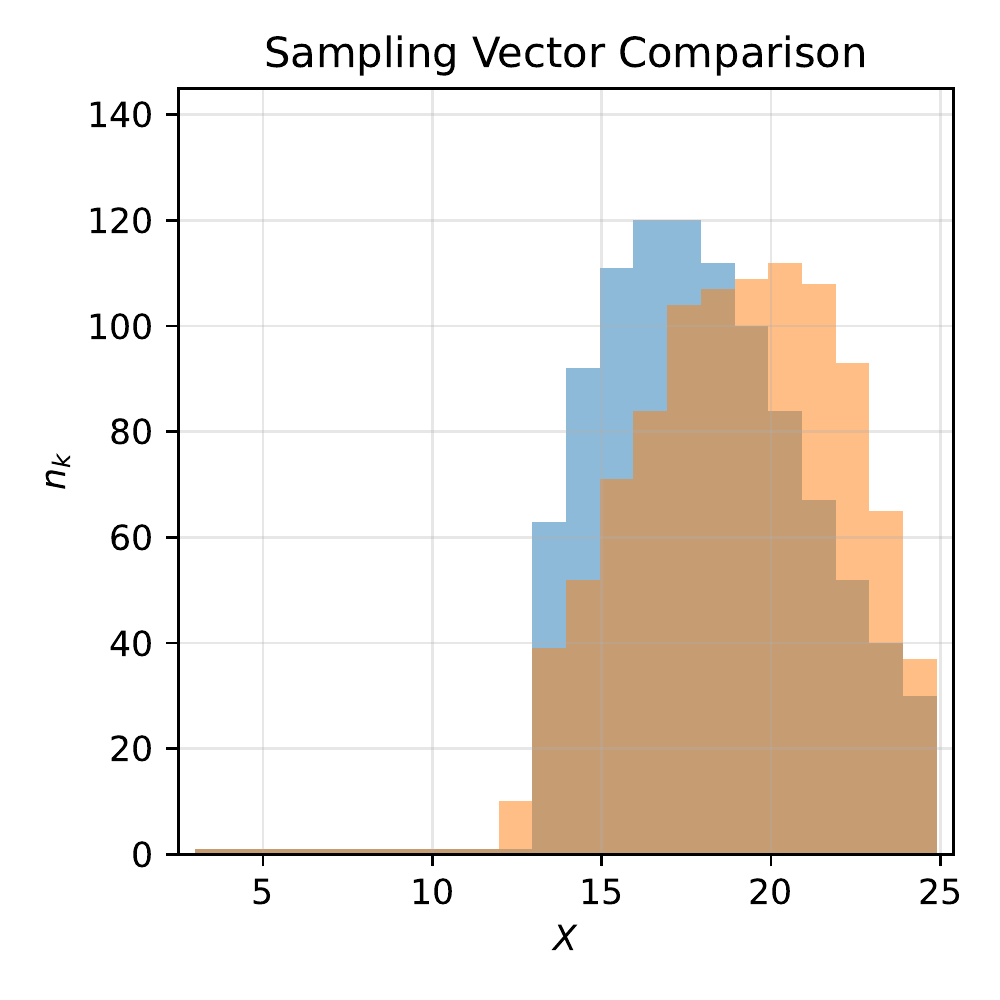}
        \caption{$1$-Wasserstein Distance}
        \label{subfig:stratification_vector_comparison_1W_cs}
    \end{subfigure}
    \\
    \begin{subfigure}[b]{0.26\textwidth}
        \centering
        \includegraphics[width=\textwidth]{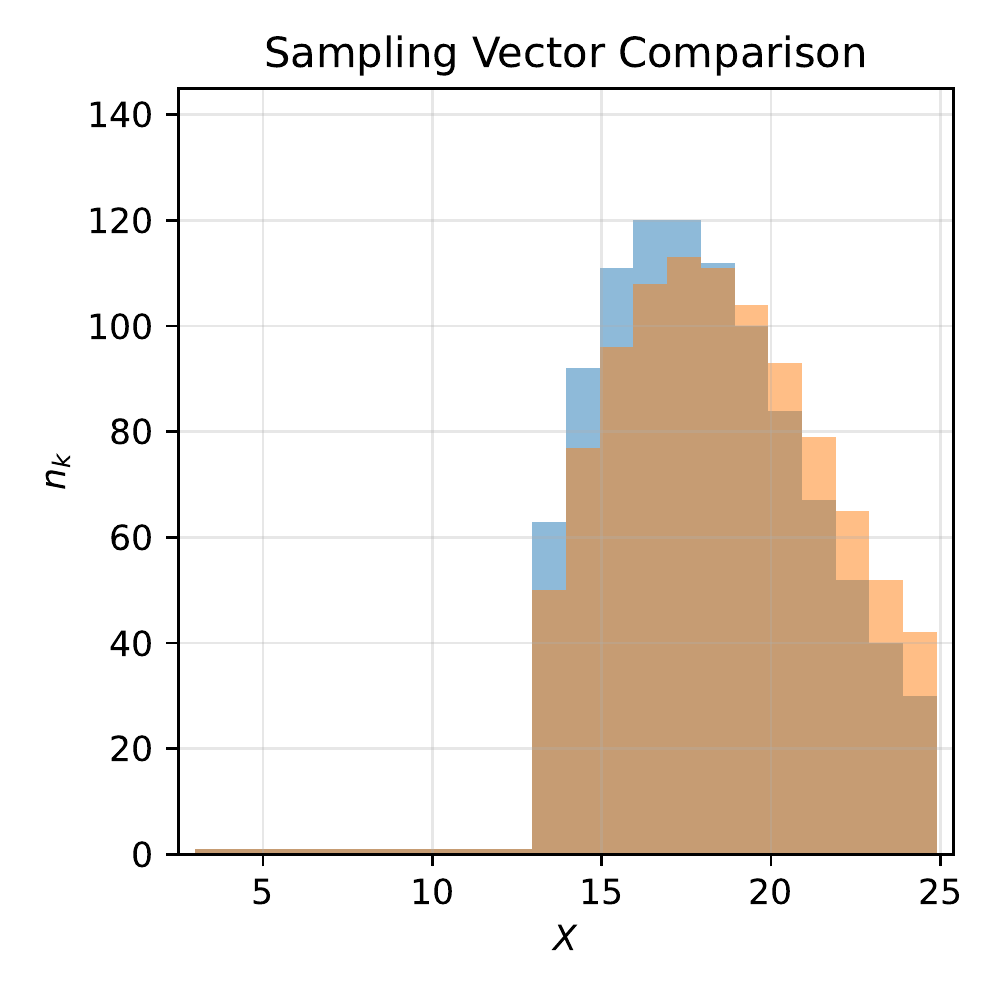}
        \caption{Parametric Family}
        \label{subfig:stratification_vector_comparison_param_cs}
    \end{subfigure}
    \begin{subfigure}[b]{0.26\textwidth}
        \centering
        \includegraphics[width=\textwidth]{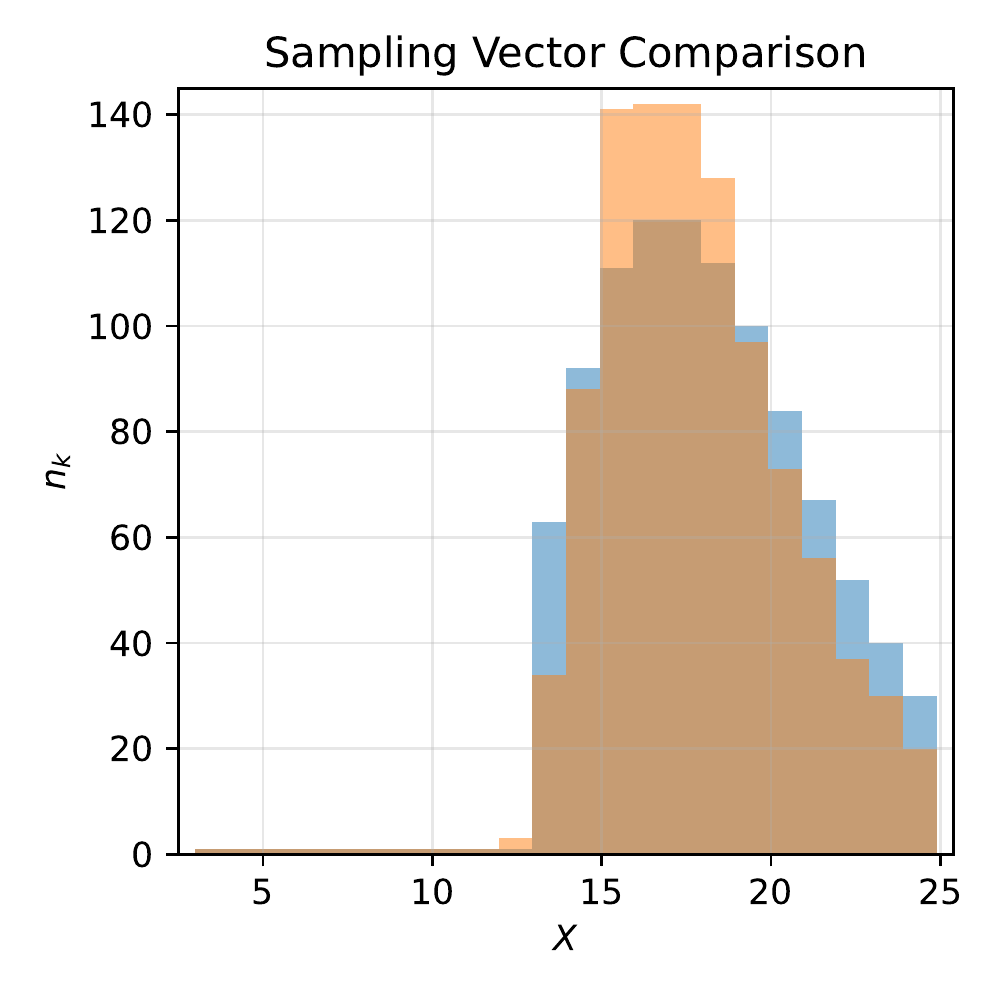}
        \caption{Distribution Moments}            \label{subfig:stratification_vector_comparison_moments_cs}
    \end{subfigure}
    \caption{Budget Allocation in Wind Turbine Case Study}
    \label{fig:stratification_vector_comparison_Strat_vs_DR_Strat_cs}
\end{figure}

In estimating the exceedance probability $\mathbb{P}(Y(X)>l)$, we set the threshold $l$ at 3.15. Because there are 220 bins, each of which is very narrow, it is not appropriate to use the bins as strata directly. Instead, we group them and take $K=22$ equally partitioned strata with $S_k = \{x_i | x_i = 3 + (k-1) + 0.1\times(i-1), \ \forall i = 1, 2,\dots, 10\}$, for $k = 1,\dots,22$ and total budget of $N_T = 1000$ simulation runs. We employ ambiguity sets in \eqref{eq:uc_l2_1d}, \eqref{eq:uc_1w_1d}, \eqref{eq:uc_rayleigh_with_shift}, and \eqref{eq:uc_moment_pmf} with the set parameters provided in Online Supplement C.2.

Figure~\ref{fig:stratification_vector_comparison_Strat_vs_DR_Strat_cs} shows the budget allocation over strata. Both approaches allocate minimal budgets for strata with $X \leq 12$ due to the rare exceedance events $\{Y(X)>l\}$ in low wind speeds. The benchmark model's sampling vector $\boldsymbol{n}^{\text{Str-M}}$ peaks around wind speed of roughly $17m/s$ where both conditional output variance and input probability are somewhat large. 

The proposed method's sampling vector $\boldsymbol{n}^{\text{DR-Str}}$ exhibits distinct patterns in different ambiguity set designs. For discrepancy-based sets, DR-strat distributes large budgets in the high wind speed region, where exceedance events are more likely to occur and the conditional output variance is higher, as shown in Figures~\ref{subfig:stratification_vector_comparison_l2_cs} and~\ref{subfig:stratification_vector_comparison_1W_cs}. Similar patterns may be seen in the budget distribution using a parametric family of ambiguity sets in Figure~\ref{subfig:stratification_vector_comparison_param_cs}, although the budgets for the right tail are smaller than those obtained using discrepancy-based ambiguity sets. The budget allocation for the moment-based ambiguity set tends to concentrate on the mid-wind speed. This is due to the fact that the conditional output variance is in unimodal form with its mode near $X=20$. The probability of the mid-wind speed region is determined to be highest when finding the inner problem's worst-case estimator variance while satisfying the moment constraints. 
    
\begin{table}[t!]
\resizebox{\textwidth}{!}{\begin{tabular}{c|c|cccc}
\toprule
& \multirow{2}{*}{Input Model \#} & \multicolumn{4}{c}{Ambiguity Set Type} \\ \cline{3-6}
&  & \multicolumn{1}{c|}{$L_2$-Norm} & \multicolumn{1}{c|}{1-Wasserstein Distance} & \multicolumn{1}{c|}{Parametric Family} & Distribution Moments \\ \hline
\multirow{2}{*}{DR-Strat} & 1  & \multicolumn{1}{c|}{$1.278\times 10^{-5}$}  & \multicolumn{1}{c|}{$3.026\times 10^{-5}$}  & \multicolumn{1}{c|}{$1.514\times 10^{-5}$}  & $1.574\times 10^{-3}$  \\ \cline{2-6} 
& 2  & \multicolumn{1}{c|}{$1.659\times 10^{-5}$}  & \multicolumn{1}{c|}{$3.901\times 10^{-5}$}  & \multicolumn{1}{c|}{$1.794\times 10^{-5}$}  & $2.344\times 10^{-3}$  \\ \hline
\multirow{2}{*}{Benchmark Model} & 1 & \multicolumn{1}{c|}{$1.747\times 10^{-5}$}  & \multicolumn{1}{c|}{$3.209\times 10^{-5}$}  & \multicolumn{1}{c|}{$1.484\times 10^{-5}$}  & $1.863\times 10^{-3}$ \\ \cline{2-6} 
& 2 & \multicolumn{1}{c|}{$2.130\times 10^{-5}$}  & \multicolumn{1}{c|}{$5.546\times 10^{-5}$}  & \multicolumn{1}{c|}{$1.836\times 10^{-5}$}  & $2.777\times 10^{-3}$ \\ \hline
\multirow{2}{*}{\begin{tabular}[c]{@{}c@{}}Ratio of the Maximum \\ Worst-Case Estimator Variances\end{tabular}} & \multirow{2}{*}{(Str-M)/(DR-Str)} & \multicolumn{1}{c|}{\multirow{2}{*}{1.284}} & \multicolumn{1}{c|}{\multirow{2}{*}{1.422}} & \multicolumn{1}{c|}{\multirow{2}{*}{1.025}} & \multirow{2}{*}{1.183} \\
& & \multicolumn{1}{c|}{} & \multicolumn{1}{c|}{}& \multicolumn{1}{c|}{}   & \\ \hline
\end{tabular}} 
\caption{Worst-Case Estimator Variances for Two Input Models in Case Study}
\label{tab:worst_case_estimator_variances_at_all_input_models_cs}
\end{table}

Table~\ref{tab:worst_case_estimator_variances_at_all_input_models_cs} compares the worst-case estimator variance in both input models. The ratio in the last row is calculated by dividing the benchmark model's maximum worst-case variance by that of the DR-strat. The DR-strat always yields a smaller worst-case variance, indicating its robustness. For the four different forms of ambiguity sets, we observe the various levels of variance reduction.  The discrepancy-based sets, followed by the moment-based set, show the greatest reduction among the four ambiguity sets. Online Supplement D.3 provides more detailed experimental results, including the pmf instances within each ambiguity set and the worst-case distributions of the inner problem.

\section{Conclusions}
\label{sec:conclusions}

This paper proposes a robust stratified sampling method to address multiple uncertain input models. We formulate an optimization problem to minimize the maximum of worst-case estimator variances among candidate distributions based on the DRO framework. We solve the resulting bi-level optimization problem using BO to obtain the robust DR-strat sampling vector, which enables the efficient reuse of simulation results.

Our numerical experiments in two settings--toy example and a case study involving a wind turbine--suggest that the proposed approach shows robust performance when the true model realization deviates from the initial belief. In comparison to the benchmark model that does not incorporate uncertainty, it obtains lower estimator variance. We also offer a thorough analysis using four different kinds of ambiguity sets and discuss how they impact the estimation outcome and under what circumstances a particular set is preferable.

Future work could investigate other variance reduction techniques, such as importance sampling and antithetic sampling, in the presence of input uncertainty. We could also explore robust simulation with multi-fidelity models. For example, we could achieve an optimal balance between estimation accuracy and simulation budget by using high-fidelity models when necessary and supplementing with cheap, low-fidelity models as needed.


\ACKNOWLEDGMENT{%
This work was supported in part by the National Research Foundation of Korea (NRF) grant funded by the Korean government (MSIT) (No. NRF-2021R1A2C1094699 and NRF-2021R1A4A1031019) and in part by U.S. National Science Foundation (CMMI-2226348 and IIS-1741166) 
}

%
%
%


\bibliographystyle{informs2014}
\bibliography{dr_strat_manuscript_JOC}


\end{document}